\documentclass[a4paper, english]{amsart}

\usepackage[utf8]{inputenc}
\usepackage[T1]{fontenc}

\usepackage{xfrac}
\usepackage{mathtools}
\usepackage{amssymb}
\usepackage{amsthm}
\usepackage{thmtools}

\usepackage{todonotes}

\usepackage{multicol}

\usepackage{tikz-cd}

\usepackage{csquotes}
\usepackage[style=alphabetic, url=false]{biblatex}

\usepackage{bbold}
\usepackage{enumitem}

\usepackage{hyperref}
\usepackage{cleveref}

\usepackage{babel}

\newcommand{\spacedequal}[0]{\quad=\quad}
\newcommand{\xspacedequal}[1]{\ensuremath{\quad\overset{#1}{=}\quad}}

\newcommand{\catname}[1]{\ensuremath{\textnormal{#1}}}
\newcommand{\bicatname}[1]{\ensuremath{\normalfont\textbf{#1}}}
\newcommand{\field}[0]{\ensuremath{\mathbb{k}}}
\newcommand{\pfield}[1]{\ensuremath{\mathbb{F}_{#1}}}

\newcommand{\catmod}[1]{\ensuremath{\catname{mod}(\field{}#1)}}

\newcommand{\catstmod}[1]{\ensuremath{\catname{stmod}(\field{}#1)}}

\newcommand{\catder}[1]{\ensuremath{\catname{D}(\field{}#1})}

\newcommand{\ncohomology}[3]{\ensuremath{H^{#1}({#2} ; {#3})}}
\newcommand{\cohomology}[2]{\ensuremath{\ncohomology{*}{#1}{#2}}}
\newcommand{\ncohomologyh}[3]{\ensuremath{\hat{H}^{#1}({#2} ; {#3})}}
\newcommand{\cohomologyh}[2]{\ensuremath{\ncohomologyh{*}{#1}{#2}}}

\newcommand{\catcat}[0]{\bicatname{Cat}}
\newcommand{\catadd}[0]{\bicatname{Add}}

\newcommand{\catgpd}{\bicatname{gpd}}
\newcommand{\catgpdf}{\bicatname{gpd}^{\textnormal{f}}}

\newcommand{\catset}{\bicatname{Set}}

\newcommand{\catfusion}[2]{\ensuremath{\mathcal{F}_{#2}(#1)}}
\newcommand{\cattransport}[2]{\ensuremath{\mathcal{T}_{#2}(#1)}}

\newcommand{\cattransportext}[2]{\ensuremath{\hat{\mathcal{T}}_{#2}(#1)}}

\newcommand{\catop}[1]{{#1}^{\normalfont\text{op}}}

\newcommand{\group}[1]{\ensuremath{#1}}
\newcommand{\groupoid}[1]{\ensuremath{#1}}
\newcommand{\mackeyfunctor}[1]{\ensuremath{\mathcal{#1}}}

\newcommand{\ffunctor}[1]{\ensuremath{\mathbb{#1}}}

\newcommand{\monad}[1]{\ensuremath{\mathsf{#1}}}

\newcommand{\ppback}[2]{\ensuremath{(#1|#2)}}
\newcommand{\ppbackl}[2]{\ppback{\underline{#1}}{#2}}
\newcommand{\ppbackr}[2]{\ppback{#1}{\underline{#2}}}
\newcommand{\ppbacknt}[2]{\overline{\ppback{#1}{#2}}}

\numberwithin{equation}{section}
\swapnumbers

\declaretheorem[style=plain,sibling=equation]{theorem}

\declaretheorem[style=plain,
sibling=equation
]{proposition}

\declaretheorem[style=definition,sibling=equation]{definition}
\declaretheorem[style=definition,sibling=equation]{example}
\declaretheorem[style=remark,sibling=equation]{remark}
\declaretheorem[style=definition, 
sibling=equation,
name=Notation,
refname={notation,notations},
Refname={Notation,Notations}
]{notation}

\DeclareMathOperator{\id}{Id}

\DeclareMathOperator*{\llim}{bilim}
\DeclareMathOperator*{\ccolim}{bicolim}
\DeclareMathOperator{\Ob}{Ob}

\usepackage{newunicodechar}

\newunicodechar{よ}{\text{\usefont{U}{min}{m}{n}\symbol{'210}}}

\DeclareFontFamily{U}{min}{}
\DeclareFontShape{U}{min}{m}{n}{<-> udmj30}{}

\title{A categorification of the Cartan-Eilenberg formula}
\author{Jun Maillard}

\address{\ \medbreak
  Univ.\ Lille, CNRS, UMR 8524 - Laboratoire Paul Painlevé, F-59000 Lille, France}
\email{jun.maillard@gmail.com}
\urladdr{http://math.univ-lille1.fr/~{}jmaillard}
\thanks{Author supported by Project ANR ChroK (ANR-16-CE40-0003) and Labex CEMPI (ANR-11-LABX-0007-01).}


\newcommand{\cone}[1]{\ensuremath{\mathcal{#1}}}

\usetikzlibrary{shapes.misc, positioning}

\tikzset{tip/.style={}}
\tikzset{dot/.style={circle, fill, minimum size=4pt, inner sep=0pt, outer sep=0pt}}
\tikzset{iden/.style={draw, circle, minimum size=4pt, inner sep=0pt, outer sep=0pt}}
\tikzset{naturaltr/.style={draw, rounded rectangle,
    minimum height=0.6cm, minimum width=1.8cm, inner sep=1mm}}

\allowdisplaybreaks

\begin{document}

\begin{abstract}
  We prove a categorification of the stable elements formula of
  Cartan and Eilenberg. Our formula expresses the derived category and the stable
  module category of a group as a bilimit of the corresponding categories for the $p$-subgroups.
\end{abstract}

\maketitle

A classical theorem by Cartan and Eilenberg~\autocite[Theorem
10.1]{cartanHomologicalAlgebra1956} presents the mod $p$ cohomology of a finite
group \group{G} as a subalgebra of the cohomology of any $p$-Sylow subgroup
\group{S} of \group{G}. The
formalism of fusion systems provides a compact formula to express this result in
terms of the fusion category $\catfusion{\group{G}}{\group{S}}$ (\cref{defFusionCat}):
\[
  \cohomology{\group{G}}{\pfield{p}} \simeq \lim_{\group{P} \in \catop{\catfusion{\group{G}}{\group{S}}}}\cohomology{\group{P}}{\pfield{p}}
\]
This formula actually holds~\autocite{parkMislinTheoremFusion2017} for any cohomological Mackey functor
\mackeyfunctor{M}:
\[
  \mackeyfunctor{M}(\group{G}) \simeq \lim_{\group{P} \in \catop{\catfusion{\group{G}}{\group{S}}}}\mackeyfunctor{M}(\group{P})
\]
We prove a categorified Cartan-Eilenberg formula (\cref{theoremMackeyAs2Lim})
for any \emph{$p$-monadic Mackey $2$-functor} $\ffunctor{M}$ (\cref{defMack,defMonadMack}):
\[
  \ffunctor{M}(G) \cong \llim_{P \in \catop{\cattransportext{G}{S}}}{\ffunctor{M}(P)}
\]
For instance, the left-hand side can stand for the categories of group representations
$\ffunctor{M}(G) = \catmod{G}$, the stable categories of group representations $\ffunctor{M}(G) = \catstmod{G}$ or the derived categories of
the group algebras $\ffunctor{M}(G) = \catder{G}$. (There are other examples,
see \cref{rmkExampleMack}.) We can let $G$
range over all finite groupoids and these mappings then describe the values on objects of
$2$-functors $\ffunctor{M} \colon \catop{(\catgpdf)} \to \catcat$ from the $2$-category of
finite groupoids, with faithful functors as $1$-morphisms, to the $2$-category
of small categories.

This categorification requires us to replace the fusion category $\catfusion{G}{S}$ by an extended
transporter category $\cattransportext{G}{S}$ (\cref{defCatTransportExt}) and
the classical limit by a pseudo bilimit (\cref{defBilimit}) taken in the
$2$-category of categories.

Through a 2-finality argument, using the criterion of \autocite{maillard2final2functors2021}, the bilimit in \cref{theoremMackeyAs2Lim} can be
reinterpreted as a descent-shaped bilimit. \Cref{theoremMackeyAs2Lim} thus
states that the functor $\ffunctor{M}$ is a $2$-sheaf, which allows us to
recover the main theorem of \autocite{balmerStacksGroupRepresentations2015}.
More details on this point can be found in~\autocite{maillardFiniteGroupRepresentations2021}.

The article is organized as follows. In \cref{secBicat}, we recall classical definitions concerning
bicategories and relevant notations. In \cref{secFus}, we give a quick summary
of the classical categorical formalism of fusion in groups and
introduce our generalization. In \cref{secTheorem}, we state and prove the main
theorem. In \cref{secApp}, we explain how to retrieve
the classical Cartan-Eilenberg formula, and similar results, from our categorified formula.

\subsection*{Acknowledgements} This work is a part of a PhD thesis under the
supervision of Ivo Dell'Ambrogio.

\section{Bicategorical notions}\label{secBicat}

\subsection{Usual definitions and conventions}

We will follow the naming conventions of \autocite{johnsonDimensionalCategories2021}
for bicategorical notions. In particular, the terms \emph{$2$-category} and \emph{$2$-functor} denote the \emph{strict}
ones; \emph{bicategory} and \emph{pseudofunctor} refer to the strong but
(possibly) non-strict variants. We will use the term \emph{$(2,1)$-category} for a
\emph{$2$-category} with only invertible \emph{$2$-morphisms}, and the term \emph{$(2,1)$-functor}
for a $2$-functor between \emph{$(2,1)$-categories}.

We recall some usual constructions and properties of $2$-categories we will use,
and introduce some notations.

\begin{notation}
  We use the symbol $\simeq$ to denote an \emph{isomorphism} between two objects (in a
  $1$-category) and the symbol $\cong$ to denote an \emph{equivalence} between two objects (in a $2$-category).
\end{notation}

\begin{definition}\label{defOpCat}
  Let $\bicatname{C}$ be a 2-category. The \emph{opposite $2$-category}
  $\catop{\bicatname{C}}$ of $\bicatname{C}$ is the $2$-category with:
  \begin{itemize}
  \item \emph{Objects:} the objects of $\bicatname{C}$
  \item \emph{Hom-categories:} $\catop{\bicatname{C}}(A, B) = \bicatname{C}(B, A)$
  \end{itemize}
\end{definition}

\begin{notation}
  We use the notation $[\bicatname{A},\bicatname{B}]$ for the $2$-category of pseudofunctors,
  pseudonatural transformations and modifications, between two $2$-categories
  $\bicatname{A}$ and $\bicatname{B}$.
\end{notation}

\begin{notation}
 Let $\bicatname{I}, \bicatname{C}$ be $2$-categories and $T$ be an object of
 $\bicatname{C}$. The constant 2-functor $\bicatname{I} \to
 \bicatname{C}$ with value $T$ is denoted by $\Delta{T}$.
\end{notation}

\begin{definition}\label{defBilimit}
 Let $\bicatname{I}, \bicatname{C}$ be $2$-categories and $\ffunctor{D} \colon
 \bicatname{I} \to \bicatname{C}$ be a $2$-functor. A \emph{(pseudo) bilimit} of
 $\ffunctor{D}$ is an object $L$ of $\bicatname{C}$ and a family of equivalences
 \[
   \Psi_{T} \colon \bicatname{C}(T, L) \cong [\bicatname{I},\bicatname{C}](\Delta{T}, \ffunctor{D})
 \]
 pseudonatural in $T$. When it exists, the bilimit of $\ffunctor{D}$ is unique
 up to equivalence and the object $L$ is noted
 $\llim_{\bicatname{I}}\ffunctor{D}$.

 We call the canonical pseudonatural transformation $\Psi_{L}(\id_L)$ the \emph{standard cone} of the
 bilimit $L$.
\end{definition}

\begin{definition}
  Let $\bicatname{C}$ be a $2$-category. Let $f \colon a \to c$ and $g \colon b
  \to c$ be two $1$-morphisms of $\bicatname{C}$ with the same target. A \emph{$2$-pullback}
  of $f$ and $g$ is a bilimit of the diagram:
  \[
    \begin{tikzcd}
     a \arrow[r, "f"] & c & b \arrow[l, "g", swap]
    \end{tikzcd}
  \]
\end{definition}

\begin{notation}\label{nota2Pullback}
 We write $\ppback{f}{g}$ for the $2$-pullback of $f$ and $g$. We will use the
 following naming scheme for the structural $1$- and $2$-morphisms:
 \[
   \begin{tikzcd}
     & c & \\
     a \arrow[ru, "f"] \arrow[rr, "\ppbacknt{f}{g}"'{inner sep=1ex}, "\sim", Rightarrow, shorten=1cm]
     & & b \arrow[lu, "g", swap] \\
     & \ppback{f}{g} \arrow[lu, "\ppbackl{f}{g}"] \arrow[ru, "\ppbackr{f}{g}", swap]
   \end{tikzcd}
 \]
\end{notation}

\begin{definition}
  Let $\bicatname{C}$ be a $(2,1)$-category. Fix a $(2,1)$-functor $\ffunctor{F}
  \colon \bicatname{I} \to \bicatname{C}$ and an object $c$ of $\bicatname{C}$. The
  \emph{slice} $\ffunctor{F} / c$ is the $(2, 1)$-category with:
  \begin{itemize}
  \item \emph{Objects:} the pairs $(i, f)$ consisting of an object $i$ of
    $\bicatname{I}$ and a morphism $f \colon \ffunctor{F}i \to c$
  \item \emph{Morphisms $(i, f) \to (i', f')$:} the pairs $(u, \mu)$ consisting
    of a morphism $u \colon i \to i'$ of $\bicatname{I}$ and a 2-isomorphism $\mu \colon f \Rightarrow
    f'\ffunctor{F}(u)$ of $\bicatname{C}$:
    \[
      \begin{tikzcd}
        \ffunctor{F}i \arrow[rd, "f"{name=n}, bend right, swap] \arrow[rr, "\ffunctor{F}u"] & &
        \ffunctor{F}i' \arrow[ld, "f'", bend left] \\
        & c \arrow[from=n, to=1-3, Rightarrow, "\mu", "\sim"', shorten=0.5cm] &
      \end{tikzcd}
    \]
  \item \emph{2-Morphisms $(u, \mu) \Rightarrow (v, \nu)$:} the 2-morphisms
    $\alpha \colon u \Rightarrow v$ of $\bicatname{I}$ satisfying:
    \[
      \begin{tikzcd}[row sep=2cm]
        \ffunctor{F}i
          \arrow[rr, "\ffunctor{F}v"{name=nv}, bend left]
          \arrow[rr, "\ffunctor{F}u"{name=nu}, bend right, swap]
          \arrow[rd, "f"{name=n}, bend right, swap] & &
        \ffunctor{F}i'
          \arrow[ld, "f'", bend left]\\
        & c
        \arrow[from=n, to=1-3, Rightarrow, "\mu", "\sim"', shorten=0.5cm, shift right=0.4cm]
        \arrow[from=nu, to=nv, Rightarrow, "\ffunctor{F}\alpha", shorten=0.2cm] &
      \end{tikzcd} =
      \begin{tikzcd}[row sep=2cm]
        \ffunctor{F}i
          \arrow[rr, "\ffunctor{F}v", bend left]
          \arrow[rd, "f"{name=n}, bend right, swap] & &
        \ffunctor{F}i'
          \arrow[ld, "f'", bend left] \\
        & c
        \arrow[from=n, to=1-3, Rightarrow, "\nu", "\sim"', shorten=0.5cm] &
      \end{tikzcd}
    \]
  \item Compositions are induced by the compositions of $\bicatname{I}$ and $\bicatname{C}$.
  \end{itemize}
  A slice $2$-category $\ffunctor{F}/c$ is endowed with a canonical forgetful $2$-functor:
  \[
    \left\{
      \begin{array}{lll}
        \ffunctor{F} / c & \to & \bicatname{I} \\
        (i, f) & \mapsto & i \\
        (u, \mu) & \mapsto & u \\
        \alpha & \mapsto & \alpha
      \end{array}
    \right.
  \]
\end{definition}

\subsection{Adjunctions and monads}

We recall the definitions of adjunctions and monads, and some of their basic
properties in $\catcat$.

\begin{definition}
 Let $C, D$ be two categories. An adjunction between $C$ and $D$ is a quadruple
 $(\ell, r, \eta, \epsilon)$ of a functor $\ell \colon C \to D$, a functor $r \colon D
 \to C$, a natural transformation $\eta \colon \id_C \Rightarrow r\ell$ and a natural
 transformation $\epsilon \colon \ell{}r \Rightarrow \id_D$, such that:
 \begin{equation}\label{eqUC1}
   \id_{\ell} = \epsilon\ell \circ \ell\eta
 \end{equation}
 \begin{equation}\label{eqUC2}
   \id_r = r\epsilon \circ \eta{}r
 \end{equation}
 We write an adjunction $\ell \dashv r$, with the natural transformations $\eta$ and
 $\epsilon$ omitted. The natural transformation $\eta$ is called the \emph{unit}
 of the adjunction and the natural transformation $\epsilon$ is called the \emph{counit}.
\end{definition}

\begin{definition}
 Let $C$ be a category. A monad on $C$ is a monoid object in the monoidal category
 $\catcat(C, C)$ of endofunctors of $C$, that is, a triple $(\monad{T}, \eta,
 \mu)$ consisting of a functor
 $\monad{T} \colon C \to C$, a natural transformation $\eta \colon \id_C \Rightarrow
 \monad{T}$ and a natural transformation $\mu \colon \monad{T} \circ \monad{T}
 \Rightarrow \monad{T}$ such that the following diagrams commute:
 \[
   \begin{tikzcd}
     \monad{T} \circ \monad{T} \circ \monad{T}
     \arrow[r, "\monad{T}\mu", Rightarrow]
     \arrow[d, "\mu\monad{T}", Rightarrow, swap] &
     \monad{T} \circ \monad{T} \arrow[d, "\mu", Rightarrow] \\
     \monad{T} \circ \monad{T} \arrow[r, "\mu", Rightarrow, swap] & \monad{T}
   \end{tikzcd} \qquad
   \begin{tikzcd}
     \monad{T} \circ \monad{T} \arrow[rd, "\mu", Rightarrow, swap] &
     \monad{T}
     \arrow[l, "\monad{T}\eta", Rightarrow, swap]
     \arrow[r, "\eta\monad{T}", Rightarrow] \arrow[d, equal] &
     \monad{T} \circ \monad{T} \arrow[ld, "\mu", Rightarrow] \\
     & \monad{T} &
   \end{tikzcd}
 \]
\end{definition}

\begin{proposition}
 Let $C, D$ be two categories and $(\ell \colon C \to D,r \colon D \to
 C,\eta,\epsilon)$ be an adjunction $\ell \dashv r$ between $C$ and $D$. Then the
 triple $(r\ell, \eta, r\epsilon\ell)$ defines a monad on $C$.
\end{proposition}

\begin{definition}\label{defEM}
 Let $C$ be a category and $\monad{T}$ be a monad on $C$. The \emph{Eilenberg-Moore
 category} $C^{\monad{T}}$ of $\monad{T}$ is the category with:
 \begin{itemize}
 \item \emph{Objects:} the pairs $(c, f)$ with $c$ an object of $C$ and $f
   \colon \monad{T}c \to c$ a morphism of $C$.
 \item \emph{Morphisms $(c,f) \to (c', f')$:} the morphisms $g \colon c \to c'$
   such that the following square commutes:
   \[
     \begin{tikzcd}
      \monad{T}c \arrow[r, "\monad{T}g"] \arrow[d, "f", swap] & \monad{T}c'
      \arrow[d, "f'"] \\
      c \arrow[r, "g", swap] & c'
     \end{tikzcd}
   \]
 \item \emph{Composition} is induced by the composition of $C$.
 \end{itemize}
 There is a canonical functor
 \[
   U \colon \left\{
     \begin{array}{lll}
       C^{\monad{T}} & \to & C \\
       (c, f) & \mapsto & c \\
       g & \mapsto & g
     \end{array}
   \right.
 \]
 and a natural transformation $\bar{\mu} \colon \monad{T}U \Rightarrow U$ with components
 \[
   \bar{\mu}_{(c,f)} = f
 \]
\end{definition}

\begin{definition}\label{defModuleMonad}
  Let $C$ be a category and $(\monad{T}, \eta, \mu)$ be a monad on $C$. A \emph{left module} on
  $\monad{T}$, or left $\monad{T}$-module, is a category $D$ endowed with a functor $V \colon D \to C$ and
  a natural transformation $\nu \colon \monad{T}V \Rightarrow V$ such that the
  following diagrams commute
  \[
    \begin{tikzcd}
      \monad{T}\monad{T}V
      \arrow[r, "\mu{}V", Rightarrow]
      \arrow[d, "\monad{T}\nu", Rightarrow, swap] &
      \monad{T}V \arrow[d, "\nu", Rightarrow] \\
      \monad{T}V \arrow[r, "\nu", Rightarrow, swap] &
      V
    \end{tikzcd} \qquad
    \begin{tikzcd}
      V \arrow[r, "\eta{}V", Rightarrow] \arrow[d, equal] &
      \monad{T}V \arrow[ld, "\nu", Rightarrow] \\
      V &
    \end{tikzcd}
  \]
\end{definition}

\begin{definition}\label{defModule2Cat}
  Let $C$ be a category and $(\monad{T}, \eta, \mu)$ be a monad on $C$.
  The \emph{$2$-category of left $\monad{T}$-modules} is the category with:
  \begin{itemize}
  \item \emph{Objects:} the left $\monad{T}$-modules.
  \item \emph{$1$-morphisms $(D, V, \nu) \to (E, W, \omega)$:}
    the pairs $(X, \chi)$ consisting of a functor $X \colon D \to E$ and a
    natural isomorphism $\chi \colon V \xRightarrow{\sim} WX$ such that:
    \[
      \chi \circ \nu = \omega{}X \circ \monad{T}\chi
    \]
  \item \emph{$2$-morphisms $(X, \chi) \Rightarrow (X', \chi')$:} the
    natural transformations $\psi : X \Rightarrow X'$ such that:
    \[
      \begin{tikzcd}[row sep=2cm]
        D \arrow[rr, bend left, "X'"{name=n1}] \arrow[rr, bend right, "X"{name=n2}, swap]
          \arrow[rd, bend right, "V"{name=nv}, swap]
        & & E \arrow[ld, bend left, "W"] \\
        & C &
        \arrow[from=n2, to=n1, "\psi", Rightarrow, shorten=0.25cm]
        \arrow[from=nv, to=1-3, "\chi", Rightarrow, swap, shorten=0.5cm, shift right=0.4cm]
      \end{tikzcd} =
      \begin{tikzcd}[row sep=2cm]
        D \arrow[rr, bend left, "X'"] \arrow[rd, bend right, "V"{name=nv}, swap]
        & & E \arrow[ld, bend left, "W"] \\
        & C &
        \arrow[from=nv, to=1-3, "\chi'", Rightarrow, swap, shorten=0.5cm]
      \end{tikzcd}
    \]
  \end{itemize}
\end{definition}

\begin{proposition}\label{propEMModule}
 Let $C$ be a category and $\monad{T}$ be a monad on $C$. The Eilenberg-Moore
 category $C^{\monad{T}}$ endowed with $U$ and $\bar{\mu}$ (\cref{defEM}) is a
 left module on $\monad{T}$ (\cref{defModuleMonad}).

 Moreover it is terminal among the left modules on $\monad{T}$ in the following sense: if $(D, V, \nu)$
 is another left module on $\monad{T}$, there is a unique functor $K \colon D
 \to C^{\monad{T}}$ such that $V = UK$ and $\nu = \bar{\mu}K$.
\end{proposition}

\begin{remark}\label{remarkEMBiterm}
 The left $\monad{T}$-module $(C^{\monad{T}}, U, \bar{\mu})$ of
 \cref{propEMModule} is pseudo biterminal in the $2$-category of left
 $\monad{T}$-modules (\cref{defModule2Cat}): for any
 other left $\monad{T}$-module $(D, V, \nu)$, there is a $1$-morphism $(K,
 \kappa) \colon (D,V,\nu) \to (C^{\monad{T}},U,\bar{\mu})$, unique up to a unique $2$-isomorphism.

\end{remark}

\subsection{String diagrams}

The proof of our main result is an extensive manipulation of string diagrams. We recall the
general ideas of string diagrams, and introduce our notations. A fully detailed
description of string diagrams can be found in \autocite[§3.7]{johnsonDimensionalCategories2021}.

String diagrams are used to depict natural transformations. They are the dual
diagrams of the ``usual'' pasting diagrams.
\begin{itemize}
\item An object $A$ is represented by a labeled surface
  \[
    \begin{tikzpicture}
     \node (a) at (0, 0) {$A$};
    \end{tikzpicture}
  \]
\item A 1-morphism $f \colon A \to B$ is represented by a labeled vertical edge
  
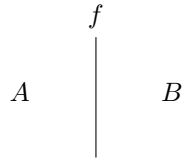
\[
    \begin{tikzpicture}
      \node[tip] (sta) at (0, 0)  {$f$};
      \node[tip] (end) at (0, -2) {};
      \node (a) at (-1, -1) {$A$};
      \node (b) at ( 1, -1) {$B$};
      \draw (sta) to (end);
    \end{tikzpicture}
  \]
  with the source on the left and the target on the right.
\item A 2-morphism $\phi \colon f \Rightarrow g \colon A \to B$ is represented
  by a labeled vertex
  
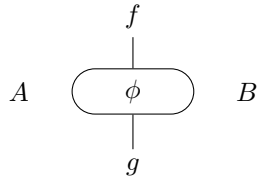
\[
    \begin{tikzpicture}
      \node[tip] (sta) at (0, 0)  {$f$};
      \node[naturaltr] (phi) at (0, -1) {$\phi$};
      \node[tip] (end) at (0, -2) {$g$};
      \node (a) at (-1.5, -1) {$A$};
      \node (b) at ( 1.5, -1) {$B$};
      \draw (sta) to (phi.north);
      \draw (phi.south) to (end);
    \end{tikzpicture}
  \]
  with the source above and the target below. The identity 1-morphisms may be omitted.
\item We will occasionally represent identity 2-morphisms by a white dot, that
  is, the string diagram
  
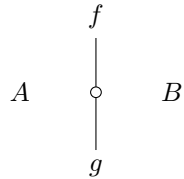
\[
    \begin{tikzpicture}
      \node[tip] (sta) at (0, 0)  {$f$};
      \node[iden] (id1) at (0, -1) {};
      \node[tip] (end) at (0, -2) {$g$};
      \node (a) at (-1, -1) {$A$};
      \node (b) at ( 1, -1) {$B$};
      \draw (sta) to (id1.north);
      \draw (id1.south) to (end);
    \end{tikzpicture}
  \]
  is the identity between the 1-morphisms $f$ and $g$.
\item When dealing with an adjunction $\ell \dashv r$, units and counits will be
  represented by black dots. For instance, the string diagram
  
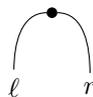
\[
    \begin{tikzpicture}
      \node[dot] (uni) at (0.5, 0) {};
      \node[tip] (end1) at (0, -1)  {$\ell$};
      \node[tip] (end2) at (1, -1) {$r$};
      \draw (uni) to[out=180, in=90] (end1);
      \draw (uni) to[out=0, in=90] (end2);
    \end{tikzpicture}
  \]
  is the unit of the adjunction $\ell \dashv r$.
\end{itemize}

\begin{remark}
 The notation for the unit and the counit of an adjunction $\ell \dashv r$ gives a simple
 interpretation of the unit-counit laws (\cref{eqUC1} and \cref{eqUC2}):
 \[
   \begin{tikzpicture}[baseline=(current bounding box.center)]
     \node[tip] (sta1) at (0, 0) {$\ell$};
     \node[tip] (end1) at (0, -1.5)  {$\ell$};
     \draw (sta1) to[out=270, in=90] (end1);
   \end{tikzpicture} \spacedequal
   \begin{tikzpicture}[baseline=(current bounding box.center)]
     \node[tip] (sta1) at (1.5, 0) {$\ell$};
     \node[dot] (uni) at (0.5, -0.5) {};
     \node[dot] (cou) at (1, -1) {};
     \node[tip] (end1) at (0, -1.5)  {$\ell$};
     \draw (sta1) to[out=270, in=0] (cou);
     \draw (uni) to[out=180, in=90] (end1);
     \draw (uni) to[out=0, in=180] (cou);
   \end{tikzpicture}
   \qquad \qquad \qquad
   \begin{tikzpicture}[baseline=(current bounding box.center)]
     \node[tip] (sta1) at (0, 0) {$r$};
     \node[tip] (end1) at (0, -1.5)  {$r$};
     \draw (sta1) to[out=270, in=90] (end1);
   \end{tikzpicture} \spacedequal
   \begin{tikzpicture}[baseline=(current bounding box.center)]
     \node[tip] (sta1) at (0, 0) {$r$};
     \node[dot] (uni) at (1, -0.5) {};
     \node[dot] (cou) at (0.5, -1) {};
     \node[tip] (end1) at (1.5, -1.5)  {$r$};
     \draw (sta1) to[out=270, in=180] (cou);
     \draw (uni) to[out=0, in=90] (end1);
     \draw (uni) to[out=180, in=0] (cou);
   \end{tikzpicture}
 \]
\end{remark}

\subsection{$1$-categories as $2$-categories}

A $1$-category can be seen as a $(2,1)$-category, with only identities as
$2$-morphisms, and a $1$-functor between $1$-categories can itself be viewed as
a $2$-functor. Hence all the constructions on $2$-categories can be applied to
$1$-categories and $1$-functors. Reciprocally there are several ways to extract
$1$-categories from bicategories.

\begin{definition}
 Let $\bicatname{C}$ be a bicategory with only invertible 2-morphisms. The \emph{truncated
 category} $\tau_1\bicatname{C}$ is the category where:
 \begin{itemize}
 \item \emph{Objects} are the objects $c$ of $\bicatname{C}$.
 \item \emph{Morphisms $c \to c'$} are the classes of $1$-morphisms $c \to c'$ of
   $\bicatname{C}$ up to $2$-isomorphism.
 \item \emph{Composition} is induced by the composition of $\bicatname{C}$.
 \end{itemize}
 There is a canonical projection pseudofunctor $\pi \colon \bicatname{C} \to \tau_1\bicatname{C}$.
\end{definition}

\begin{proposition}
 Let $\bicatname{C}$ be a bicategory with only invertible $2$-morphisms and
 $\bicatname{D}$ be a $1$-category seen as a $2$-category. Let $\ffunctor{F}
 \colon \bicatname{C} \to \bicatname{D}$ be a pseudofunctor. Then there is a unique
 $1$-functor $\tau_1\ffunctor{F} \colon \tau_1\bicatname{C} \to \bicatname{D}$
 factoring $\ffunctor{F}$ through $\pi$:
 \[
   \begin{tikzcd}
     \bicatname{C} \arrow[rd, "\ffunctor{F}"] \arrow[d, "\pi", swap] & \\
     \tau_1\bicatname{C} \arrow[r, "\tau_1\ffunctor{F}", swap] & \bicatname{D}
   \end{tikzcd}
 \]
\end{proposition}

\begin{definition}
 Let $\bicatname{C}$ be a $(2,1)$-category. The \emph{underlying category}
 $\bicatname{C}^{(1)}$ is the category where:
 \begin{itemize}
 \item \emph{Objects} are the objects $c$ of $\bicatname{C}$.
 \item \emph{Morphisms $c \to c'$} are the $1$-morphisms of
   $\bicatname{C}$.
 \item \emph{Composition} is induced by the composition of $\bicatname{C}$.
 \end{itemize}
 There is a canonical inclusion $2$-functor $\iota \colon \bicatname{C}^{(1)} \to \bicatname{C}$.
\end{definition}

\begin{proposition}\label{propBilim1Cat}
 Let $\bicatname{C}$ be a $(2,1)$-category and
 $\catname{D}$ be a $1$-category. Let $\ffunctor{F}
 \colon \bicatname{C} \to \catname{D}$ be a $2$-functor. Then the bilimit of
 $\ffunctor{F}$, if it exists, can be expressed as either of the following two
 ordinary limits in $\catname{D}$:
 \[
   \llim_{\bicatname{C}}\ffunctor{F} \simeq
   \lim_{\tau_1\bicatname{C}}\tau_1\ffunctor{F} \simeq
   \lim_{\bicatname{C}^{(1)}}{\ffunctor{F} \circ \iota}
 \]
\end{proposition}

\section{A 2-categorical framework for fusion theory} \label{secFus}

\subsection{The classical Cartan-Eilenberg formula}

\begin{definition}\label{defFusionCat}
  Let $G$ be a finite group and $S$ be a $p$-Sylow subgroup of $G$. The \emph{fusion system}
  $\catfusion{G}{S}$ of $G$ is the category with:
  \begin{itemize}
  \item \emph{Objects:} the subgroups $P$ of $S$.
  \item \emph{Morphisms $P \to Q$:} the group morphisms $P \to Q$ induced by the
    conjugation action of $G$:
    \[
      \catfusion{G}{S}(P, Q) = \{ c_g : p \mapsto gpg^{-1}
        \: | \: g \in G, gPg^{-1} \subset Q \}.
      \]
    \item \emph{Composition} is the composition of group morphisms.
  \end{itemize}
\end{definition}

The fusion system $\catfusion{G}{S}$ of a group $G$ is canonically endowed with
a forgetful functor to the category $\text{gp}$ of finite groups:
\[
  U \colon \left\{
    \begin{array}{lll}
      \catfusion{G}{S} & \to & \text{gp} \\
      P & \mapsto & P \\
      f & \mapsto & f
    \end{array}
  \right.
\]

The Cartan-Eilenberg formula~\autocite[Theorem 10.1]{cartanHomologicalAlgebra1956} expresses the cohomology (with trivial
coefficients) of a group $G$ as a limit over the fusion system
$\catfusion{G}{S}$ of $G$ in the category $\pfield{p}\text{Alg}^{\text{gr}}$ of
graded $\pfield{p}$-algebras:
\begin{proposition}
  For any finite group $G$ and $p$-Sylow $S$ of $G$:
  \[
    \cohomology{\group{G}}{\pfield{p}} \simeq
    \lim_{\catfusion{\group{G}}{\group{S}}^{op}}\cohomology{-}{\pfield{p}}
    \circ U
  \]
\end{proposition}

\begin{definition}\label{defCatTransport}
 Let $G$ be a finite group and $S$ be a $p$-Sylow of $G$. The \emph{transporter
   category} $\cattransport{G}{S}$ is the category with:
 \begin{itemize}
 \item \emph{Objects:} the subgroups $P$ of $S$.
 \item \emph{Morphisms $P \to Q$:} the elements $g$ of $G$ such that $gPg^{-1}
   \subset Q$:
   \[
     \cattransport{G}{S}(P, Q) = \{ g \in G \: | \: gPg^{-1} \subset Q \}.
   \]
 \item \emph{Composition} is given by the multiplication of $G$.
 \end{itemize}
\end{definition}

The transporter category $\cattransport{G}{S}$ of a group $G$ is also
canonically endowed with a forgetful functor to the category $\text{gp}$ of
finite groups:
\[
  U \colon \left\{
    \begin{array}{lll}
      \catfusion{G}{S} & \to & \text{gp} \\
      P & \mapsto & P \\
      g & \mapsto & c_g : p \mapsto gpg^{-1}
    \end{array}
  \right.
\]

This forgetful functor $U \colon \cattransport{G}{S} \to \text{gp}$ factors
through $\catfusion{G}{S}$. The comparison functor $\pi \colon \cattransport{G}{S} \to
\catfusion{G}{S}$ is the identity on objects and full.
\[
  \begin{tikzcd}
    \cattransport{G}{S} \arrow[rr, "\pi"] \arrow[rd, "U", swap] & & \catfusion{G}{S}
    \arrow[ld, "U"] \\
    & \text{gp} &
  \end{tikzcd}
\]

Note that the functor $\pi$ is a final functor, in the sense of \autocite[§IX.3]{maclanesaundersCategoriesWorkingMathematician1971}, since it is the identity on objects
and full. Hence, there is a canonical isomorphism:
\[
  \lim_{\catfusion{\group{G}}{\group{S}}^{op}}\cohomology{-}{\pfield{p}}
  \circ U  \simeq
  \lim_{\cattransport{\group{G}}{\group{S}}^{op}}\cohomology{-}{\pfield{p}} \circ U
\]

In particular, the Cartan-Eilenberg formula can equivalently be stated as:
\[
  \cohomology{\group{G}}{\pfield{p}} \simeq
  \lim_{\cattransport{\group{G}}{\group{S}}^{op}}\cohomology{-}{\pfield{p}}
  \circ U
\]

\subsection{The extended transporter (2,1)-category of a group}

We present a $2$-categorification $\cattransportext{G}{S}$ of the transporter
category. This $2$-category will have a role similar to the one the fusion
category $\catfusion{G}{S}$ and the transporter category $\cattransport{G}{S}$
have in the classical Cartan-Eilenberg formula.

Mostly for convenience, rather than restricting our attention to the category of
finite groups $\text{gp}$, we should consider the $(2,1)$-category of finite
groupoids $\catgpd$. It inherits all its structure of $2$-category from the
usual structure of $2$-category on $\catcat$, the category of small categories,
functors and natural transformations. The finite groups are viewed in $\catgpd$
as the groupoids with exactly one object. We should also note that the
$2$-morphisms of $\catgpd$ are relevant from the point of view of the fusion in
groups. Indeed, given two parallel morphisms $\phi,\psi \colon H \to G$ between
groups, a $2$-morphism $\phi \Rightarrow \psi$ is precisely an element $g$ of $G$
such that:
\[
  \forall h \in H, \phi(h) = g \psi(h) g^{-1}
\]
The 2-full subcategory of $\catgpd$ with faithful
functors as $1$-morphisms is denoted by $\catgpdf$.

\begin{definition}\label{defCatTransportExt}
  Let $G$ be a group and $S$ be a $p$-Sylow of $G$. Denote by $i \colon S \to G$
  the inclusion, seen as a $1$-morphism of $\catgpdf$. The \emph{extended
    transporter category} $\cattransportext{G}{S}$ is the $(2,1)$-category with:
  \begin{itemize}
  \item \emph{Objects:} the pairs $(\groupoid{P}, j_{\groupoid{P}} \colon \groupoid{P} \to S)$ of
    a finite groupoid  $\groupoid{P}$ and a faithful functor $j_{\groupoid{P}}
    \colon \groupoid{P} \to S$.
  \item \emph{1-Morphisms $(\groupoid{P}, j_{\groupoid{P}}) \to (\groupoid{Q},
      j_{\groupoid{Q}})$:} the pairs $(a, \alpha)$ of a faithful functor $a \colon
    \groupoid{P} \to \groupoid{Q}$ and $2$-morphisms $\alpha \colon
    ij_{\groupoid{Q}}a \Rightarrow ij_{\groupoid{P}}$ in $\catgpdf$.
  \item \emph{2-Morphisms $(a, \alpha) \Rightarrow (b, \beta)$:} the $2$-morphisms $\phi
    \colon \alpha \Rightarrow \beta$
    in $\catgpdf$, such that:
    \[
      \begin{tikzcd}
        & \group{G} & \\
        \groupoid{P} \arrow[ru, "ij_{\groupoid{P}}"]
        \arrow[rr, "b"{name=b}, swap, bend left] \arrow[rr,
        "a"{name=a}, swap, bend right]
        &
        {}
        \arrow[u, "\phi", Rightarrow, from=a, to=b, shorten=1mm]
        \arrow[u, "\beta",Rightarrow, from=b, to=u, shorten=1mm]
        & \groupoid{Q} \arrow[lu, "ij_{\groupoid{Q}}", swap]
      \end{tikzcd}
      =
      \begin{tikzcd}[row sep=large]
        & \group{G} & \\
        \groupoid{P} \arrow[ru, "ij_{\groupoid{P}}"] \arrow[rr, "a"{name=a}, swap, bend right]
        &
        {}
        \arrow[u, "\alpha", Rightarrow, from=a, shorten=2mm]
        & \groupoid{Q} \arrow[lu, "ij_{\groupoid{Q}}", swap]
      \end{tikzcd}
    \]
  \item \emph{Compositions} are induced by the obvious pasting of diagrams in $\catgpdf$.
  \end{itemize}
  The extended transporter category is endowed with a forgetful $2$-functor:
  \[
    \ffunctor{U} \colon \left\{
      \begin{array}{lll}
        \cattransportext{G}{S} & \to & \catgpdf \\
        (\groupoid{P}, j_{\groupoid{P}}) & \mapsto & \groupoid{P} \\
        (a , \alpha) & \mapsto & a \\
        \phi & \mapsto & \phi
      \end{array}
    \right.
  \]
\end{definition}

\begin{remark}
 Since $i \colon S \to G$ is a faithful functor between 1-object groupoids,
 the factorization of a functor $i_{\groupoid{P}} \colon \groupoid{P} \to G$
 through $i$, if it exists, is unique. Hence the objects of the transporter $2$-category
 $\cattransportext{G}{S}$ could equivalently be described as pairs
 $(\groupoid{P}, i_{\groupoid{P}} \colon \groupoid{P} \to G)$ such that
   the faithful functor $i_{\groupoid{P}}$ factors through $i$:
   \[
     i_{\groupoid{P}} = ij_{\groupoid{P}} \text{ for some } j_{\groupoid{P}}
       \colon \groupoid{P} \to S
   \]
\end{remark}

\begin{remark}
 The underlying $1$-category $\cattransportext{G}{S}^{(1)}$ is quite similar to
 the classical $\cattransport{G}{S}$. Actually, the main difference is the
 addition of finite coproducts. In \cref{secApp}, we choose not to
 distinguish them. Since we are working with product-preserving functors (over
 the opposite categories), this is not really a problem.
\end{remark}

\begin{remark}
  Biequivalently, the extended transporter category $\cattransportext{G}{S}$ can
  be defined as follows, using generic constructions of bicategories. Consider the slice (2,1)-category
  $\catgpdf / \group{G}$. Then we can define the extended transporter
  category \cattransportext{G}{S} as the 2-category of subobjects of $(S,i)$ in
  $\catgpdf / \group{G}$, that is, the full, 2-full subcategory of
  $\catgpdf / \group{G}$ over objects $(\groupoid{P}, {i_{\groupoid{P}}
    \colon \groupoid{P} \to \group{G}})$ such that $i_{\groupoid{P}}$ factors (up to
  some isomorphism) through $i \colon \group{S} \to \group{G}$:
  \[
    \begin{tikzcd}
      \groupoid{P} \arrow[rr] \arrow[rd, "i_{\groupoid{P}}", swap] & {} \arrow[d, phantom, "\exists\cong"] & \group{S} \arrow[ld, "i"] \\
      {} & \group{G}
    \end{tikzcd}
  \]
  In this construction, the forgetful $2$-functor
  \[
    \ffunctor{U} \colon \cattransportext{G}{S} \to \catgpdf
  \]
  appears as the restriction of the forgetful $2$-functor $\catgpdf /
  \group{G} \to \catgpdf$.
\end{remark}

\begin{remark}
  Note that there is a canonical way to recover the group $\group{G}$ from the
  $2$-functor $\ffunctor{U}$ as a pseudo bicolimit:
  \[
    \group{G} \cong \ccolim_{\cattransportext{G}{S}}\ffunctor{U}
  \]
  Indeed there is a canonical way to endow $G$ with a cocone structure over $\ffunctor{U}$, hence to define a comparison morphism from the bilimit to $G$.  
  Then, by noticing that the endomorphisms of the trivial group $(\{e\}, \{e\} \to S)$ in $\cattransportext{G}{S}$ is precisely the group $G$, one can construct another morphism from $G$ to the bilimit, which is a pseudoinverse of the canonical comparison morphism. A detailed proof can be found in~\autocite[Proposition 4.2.9]{maillardFiniteGroupRepresentations2021}.
\end{remark}

\section{Group representations as bilimits} \label{secTheorem}

\subsection{Group representations as Mackey 2-functors}

In this section, $\field$ is a commutative $\mathbb{Z}_{(p)}$-algebra (that is, every
prime integer different from $p$ is invertible in $\field$), for instance a
field of characteristic $p$.

\begin{notation}
  Let $G$ be a finite group.
  \begin{itemize}
  \item The category $\catmod{G}$ is the category of $\field$-linear
    representations of finite dimension of the group $G$.
  \item The category $\catstmod{G}$ is the stable category of $\field$-linear
    representations of finite dimension of the group $G$, for $\field$ a
    field of characteristic $p$.
  \item The category $\catder{G}$ is the (bounded) derived category of the
    finite dimensional $\field{}G$-modules.
  \end{itemize}
\end{notation}

The three mappings $G \mapsto \catmod{G}$, $G \mapsto \catstmod{G}$ and $G
\mapsto \catder{G}$ extend to $2$-functors $\ffunctor{M} \colon \catop{(\catgpdf)} \to \catadd$ from
the opposite of the $2$-category of finite groupoids (with faithful $1$-morphisms), to the $2$-category of
small additive categories. The image $\ffunctor{M}(f)$ of a group morphism $f \colon H \to G$,
noted $f^*$, is the \emph{restriction} along $f$; similarly, given a $2$-morphism
$\alpha \colon f \Rightarrow g$ in $\catgpdf$,  $\alpha^* \colon f^*
\Rightarrow g^*$ denotes the image $\ffunctor{M}(\alpha)$. We can see that these three $2$-functors are
\emph{cohomological Mackey $2$-functors with values in $\field{}$-linear idempotent-complete
  categories}, as we now explain.

\begin{definition} \label{defMack}
  A \emph{Mackey $2$-functor} is a $2$-functor $\ffunctor{M} \colon \catop{(\catgpdf)} \to
  \catadd$ satisfying the following four axioms:
  \begin{enumerate}[label=(Mack \arabic*)]
  \item \label{defMack1} \emph{Additivity:} If $i_G \colon G \to G \sqcup H$ and $i_H \colon H \to G \sqcup H$ are
    the canonical inclusions into a coproduct of groupoids, then the induced functor
    \[
      (i_G^*, i_H^*) \colon \ffunctor{M}(G \sqcup H) \to
      \ffunctor{M}(G) \times \ffunctor{M}(H)
    \]
    is an equivalence.
  \item \label{defMack2} If $f \colon H \to G$ is a $1$-morphism in $\catgpdf$, its image
    $f^*$ has both a left adjoint $f_!$ and a right adjoint $f_*$:
    \[
      f_! \dashv f^* \dashv f_*
    \]
  \item \label{defMack3} \emph{Beck-Chevalley property:} For any 2-pullback square in $\catgpdf$ (see \cref{nota2Pullback})
    \[
      \begin{tikzcd}
        K \arrow[r, "f"] \arrow[rd, Rightarrow, shorten=1mm, "\lambda"] & G \\
        (f|g) \arrow[u, "\ppbackl{f}{g}"] \arrow[r, "\ppbackr{f}{g}", swap] & L
        \arrow[u, "g", swap]
      \end{tikzcd}
    \]
    the mates $\lambda_! \colon \ppbackr{f}{g}_!\ppbackl{f}{g}^* \Rightarrow g^*f_!$ of
    $\lambda^*$ and $(\lambda^{-1})_* \colon g^*f_* \Rightarrow \ppbackr{f}{g}_*\ppbackl{f}{g}^*$ of $(\lambda^{-1})^*$ are invertible, where:
    \begin{equation}\label{eqDefMate}
      \lambda_! =
      \begin{tikzpicture}[baseline=(current bounding box.center)]
        \node[tip] (sta1) at (2, 0) {$\ppbackl{f}{g}^*$};
        \node[tip] (sta2) at (3, 0) {$\ppbackr{f}{g}_!$};
        \node[tip] (end1) at (0, -4) {$f_!$};
        \node[tip] (end2) at (1, -4) {$g^*$};
        \node[dot] (uni1) at (0.5, -1) {};
        \node[naturaltr] (lam) at (1.5, -2) {$\lambda^*$};
        \node[dot] (cou1) at (2.5, -3) {};
        \draw (sta1) to[out=270, in=90] (lam.north east);
        \draw (sta2) to[out=270, in=0] (cou1.center);
        \draw (uni1.center) to[out=180, in=90] (end1);
        \draw (uni1.center) to[out=0, in=90] (lam.north west);
        \draw (lam.south west) to[out=270, in=90] (end2);
        \draw (lam.south east) to[out=270, in=180] (cou1.center);
      \end{tikzpicture}
      \qquad
      \lambda_* =
      \begin{tikzpicture}[baseline=(current bounding box.center)]
        \node[tip] (sta1) at (0, 0) {$f_*$};
        \node[tip] (sta2) at (1, 0) {$g^*$};
        \node[tip] (end1) at (2, -4) {$\ppbackl{f}{g}^*$};
        \node[tip] (end2) at (3, -4) {$\ppbackr{f}{g}_*$};
        \node[dot] (uni1) at (2.5, -1) {};
        \node[naturaltr] (lam) at (1.5, -2) {$(\lambda^{-1})^*$};
        \node[dot] (cou1) at (0.5, -3) {};
        \draw (sta1) to[out=270, in=180] (cou1.center);
        \draw (sta2) to[out=270, in=90] (lam.north west);
        \draw (uni1.center) to[out=180, in=90] (lam.north east);
        \draw (uni1.center) to[out=0, in=90] (end2);
        \draw (lam.south west) to[out=270, in=0] (cou1.center);
        \draw (lam.south east) to[out=270, in=90] (end1);
      \end{tikzpicture}
    \end{equation}
  \item \label{defMack4} \emph{Ambidexterity:} For any $1$-morphism $f
    \colon H \to G$ in $\catgpdf$, the left and right adjoints of $f^*$ are
    isomorphic:
    \[
      f_! \simeq f_*
    \]
  \end{enumerate}
\end{definition}

\begin{notation}\label{notaBC}
 A $2$-functor $\ffunctor{M} \colon \catop{(\catgpdf)} \to \catcat$ which only
 satisfies the axioms \ref{defMack2} and \ref{defMack3} is said to have the
 \emph{Beck-Chevalley property}.

 A $2$-functor $\ffunctor{M} \colon \catop{(\catgpdf)} \to \catcat$ which only
 satisfies one half of the axioms \ref{defMack2} and \ref{defMack3} is said to have the
 \emph{left Beck-Chevalley property} or \emph{right Beck-Chevalley property}, accordingly.
\end{notation}

\begin{definition} \label{defCohMack}
 A \emph{cohomological} Mackey $2$-functor is a Mackey $2$-functor
 $\ffunctor{M}$ such that for any inclusion of groups $i \colon H \to G$ the
 composite natural transformation
 \[
   \begin{tikzcd}
    \id_{\ffunctor{M}(G)} \arrow[r, "\eta", Rightarrow] & i_*i^* \simeq i_!i^*  \arrow[r,
    "\epsilon", Rightarrow] & \id_{\ffunctor{M}(G)}
   \end{tikzcd}
 \]
 is equal to $[G \: : \: H]\sigma$, for some $2$-isomorphism $\sigma$. In
 particular, if $\ffunctor{M}$ takes values in $\field$-linear categories, the
 composite natural transformation $\epsilon\eta$ above is invertible whenever
 $[G \: : \: H]$ is coprime to $p$, since $\field$ is a $\mathbb{Z}_{(p)}$-algebra.
\end{definition}

The following comparison theorem between $\ffunctor{M}(G)$ and the associated
Eilenberg-Moore category appears in
\autocite{balmerCohomologicalMackey2functors2021}; we reproduce here its short
proof.
\begin{theorem}\label{theoremCohMonadic}
  Let $\ffunctor{M}$ be a cohomological Mackey $2$-functor (\cref{defMack} and
  \cref{defCohMack}) with values in $\field$-linear idempotent-complete
  categories.

  Then, for any finite group $G$ and subgroup $i \colon H \to G$ with $[G \: :
  \: H]$ coprime to $p$, the adjunction $i_{!}
  \dashv i^{*}$ is monadic. That is, the canonical comparison functor
  \[
    \ffunctor{M}(G) \to \ffunctor{M}(H)^{\monad{T}}
  \]
  is an equivalence, where $\monad{T} = i^*i_!$ is the monad induced by $i$.
\end{theorem}
\begin{proof}
  Recall that $\epsilon \colon i_!i^* \Rightarrow \id_{\ffunctor{M}(G)}$ is the
  counit of the adjunction $i_! \dashv i^*$ and
  ${\eta \colon \id_{\ffunctor{M}(G)} \Rightarrow i_*i^*}$ is the unit of the
  adjunction $i^* \dashv i_*$.

  Since $\ffunctor{M}$ is cohomological and  $[G \: : \: H]$ is coprime to $p$, the natural transformation
  $\epsilon\eta$ is invertible. Hence $\epsilon$ admits a section, and
  the adjunction $i_! \dashv i^*$ is monadic, by
  \autocite[Lemma 2.10]{balmerStacksGroupRepresentations2015}.
\end{proof}

This naturally leads to the following definition:
\begin{definition}\label{defMonadMack}
  A \emph{$p$-monadic} Mackey $2$-functor is a Mackey $2$-functor $\ffunctor{M}$ such
  that, for any inclusion of groups $i \colon H \to G$ with index
  coprime to $p$, the canonical comparison functor
  \[
    \ffunctor{M}(G) \to \ffunctor{M}(H)^{\monad{T}}
  \]
  is an equivalence, where $\monad{T} = i^*i_!$ is the monad induced by $i$.
\end{definition}

\begin{remark} \label{rmkExampleMack}
  By \cref{theoremCohMonadic}, any cohomological Mackey $2$-functor with values in
  $\field$-linear idempotent-complete categories is $p$-monadic. This includes
  the $2$-functors $\catmod{-}$, $\catstmod{-}$ and $\catder{-}$; other examples
  can be found in \autocite[Chapter~4]{balmerMackey2FunctorsMackey2020}.

  For instance, the mapping $G \mapsto \textnormal{coMack}_{\field}(G)$ associating to a
  group $G$ the category of cohomological $G$-local Mackey 1-functors extends to
  a cohomological Mackey $2$-functor. This is a corollary of
  \autocite[Proposition~7.3.2]{balmerMackey2FunctorsMackey2020} applied to the
  cohomological Mackey $2$-functor $S(G) = \textnormal{perm}_{\field{}}(G)$ of
  permutations $\field{}G$-modules and the Yoshida theorem\autocite{webbGuideMackeyFunctors2000}, characterizing the category of cohomological
  Mackey 1-functors as
  \[
    \textnormal{coMack}_{\field}(G) \cong \textnormal{Fun}_+(\textnormal{perm}_{\field{}}(G), \textnormal{Ab}).
  \]
\end{remark}

\subsection{The main theorem}

\begin{theorem}\label{theoremMackeyAs2Lim}
 Let $\ffunctor{M} \colon \catop{(\catgpdf)} \to \catadd$ be a p-monadic Mackey
 2-functor (\cref{defMonadMack}). Then, for any finite group $G$ and any $p$-Sylow $S$
 of $G$, there is an equivalence
 \[
   \ffunctor{M}(G) \cong \llim_{\catop{\cattransportext{G}{S}}}{\ffunctor{M}
     \circ \ffunctor{U}}
 \]
 where $\cattransportext{G}{S}$ is the extended transporter $2$-category of $G$ (\cref{defCatTransportExt}).
\end{theorem}

Using the $p$-monadicity, \cref{theoremMackeyAs2Lim} is a corollary of the following result:
\begin{theorem}\label{theoremBR}
  Let $\ffunctor{M} \colon \catop{(\catgpdf)} \to \catcat$ be a $2$-functor with
  the left Beck-Chevalley property (\cref{notaBC}).

  For any finite group $G$ with a $p$-Sylow subgroup $i \colon S \to G$, there is
  an equivalence
  \[
    \ffunctor{M}(S)^{\monad{T}} \cong \llim_{\catop{\cattransportext{G}{S}}}{\ffunctor{M}
      \circ \ffunctor{U}}
  \]
  where $\monad{T} = i^*i_!$ is the monad induced by $i$.
\end{theorem}

\begin{remark}
  The equivalence of \cref{theoremBR} is reminiscent of the Bénabou-Roubaud
  theorem~\autocite{benabouMonadesDescente1970}. This is actually an instance of
  a generalization to $2$-categories; a more precise statement is in~\autocite[Theorem 3.2.1]{maillardFiniteGroupRepresentations2021}. Nonetheless, the essential arguments of the proof
  are already present in this article.
\end{remark}

We will now construct the equivalence of \cref{theoremBR}. Let $L$ be the bilimit:
\[
 L = \llim_{\catop{\cattransportext{G}{S}}}{\ffunctor{M} \circ \ffunctor{U}}
\]
Its standard cone (\cref{defBilimit}) is denoted by $\cone{L}_{\bullet}$.
In the rest of this section, we will:
\begin{itemize}
\item construct a functor $A \colon L \to \ffunctor{M}(S)^{\monad{T}}$,
\item construct a functor $B \colon \ffunctor{M}(S)^{\monad{T}} \to L$,
  and
\item show that the functors $A$ and $B$ are mutual pseudoinverses.
\end{itemize}
This will prove \cref{theoremBR}.
\subsection{Construction of a comparison functor $L \to \ffunctor{M}(S)^{\monad{T}}$}

Recall that $\ffunctor{M}(S)^{\monad{T}}$ has a canonical structure of left module
on $\monad{T}$ consisting of the forgetful functor $U
: \ffunctor{M}(S)^{\monad{T}} \to \ffunctor{M}(S)$ and a natural
transformation $\bar{\mu} \colon \monad{T}U \Rightarrow U$.
Moreover, it is terminal among the left modules on $\monad{T}$
(\cref{propEMModule}).

One can define a structure of left $\monad{T}$-module for the bilimit $L$ as
follows. There is an obvious
functor $L \to \ffunctor{M}(\group{S})$, namely the component
$\cone{L}_{\group{S}}$ of the standard cone. The
definition of the action $\nu$ use the invertibility of the mate $\lambda_!$
(see \ref{defMack3}):
\begin{equation}\label{defNu}
  \begin{tikzpicture}[baseline=(current bounding box.center)]
    \node[tip] (sta1) at (0,0) {$\cone{L}_{\group{S}}$};
    \node[tip] (sta2) at (1,0) {$\monad{T}$};
    \node[tip] (end) at (0.5, -2) {$\cone{L}_{\group{S}}$};
    \node[naturaltr] (nu) at (0.5, -1) {$\nu$};
    \draw (sta1) to[out=270, in=90] (nu.north west);
    \draw (sta2) to[out=270, in=90]  (nu.north east);
    \draw (nu.south) to (end);
  \end{tikzpicture}
  \spacedequal
  \begin{tikzpicture}[baseline=(current bounding box.center)]
    \node[tip] (sta1) at (0,0) {$\cone{L}_{\group{S}}$};
    \node[tip] (sta2) at (1,0) {$i_!$};
    \node[tip] (sta3) at (2,0) {$i^*$};
    \node[naturaltr] (lam1) at (1.5, -1) {$\lambda_!^{-1}$};
    \node[naturaltr] (lam2) at (0.5,-3) {$\cone{L}_\lambda$};
    \node[dot] (cou2) at (1.5,-4) {};
    \node[tip] (end1) at (0,-5) {$\cone{L}_{\group{S}}$};
    \draw (sta1)            to[out=270, in=90] (lam2.north west);
    \draw (sta2)            to[out=270, in=90] (lam1.north west);
    \draw (sta3)            to[out=270, in=90] (lam1.north east);
    \draw (lam1.south west) to[out=270, in=90] (lam2.north east);
    \draw (lam1.south east) to[out=270, in=0]  node[right]{$\ppbackr{i}{i}_!$} (cou2.center);
    \draw (lam2.south west)        to[out=270, in=90]  (end1);
    \draw (lam2.south east)        to[out=270, in=180] node[below left]{$\ppbackr{i}{i}^*$}  (cou2.center);
\end{tikzpicture}
\end{equation}
 where $\lambda = \ppbacknt{i}{i}$ is defined by the following 2-pullback square in $\catgpdf$
\[
  \begin{tikzcd}
    & G & \\
    S \arrow[ru, "i"]
    \arrow[rr, "\lambda", Rightarrow, shorten=10mm] & & S \arrow[lu, "i", swap] \\
    & \ppback{i}{i} \arrow[lu, "\ppbackl{i}{i}"] \arrow[ru, "\ppbackr{i}{i}", swap] &
  \end{tikzcd}
\]
and $\cone{L}_{\lambda}$ is the composite $2$-morphism
\[
  \begin{tikzcd}
    \ppbackl{i}{i}^*\cone{L}_{\group{S}}
    \arrow[r, Rightarrow, "\cone{L}_{(\ppbackl{i}{i}, \id)}"]
    & \cone{L}_{\ppback{i}{i}}
    \arrow[r,Rightarrow, "\cone{L}_{(\id, \lambda)}"]
    & \cone{L}_{\ppback{i}{i}}
    \arrow[r, Rightarrow, "\cone{L}_{(\ppbackr{i}{i}, \id)}^{-1}"]
    & \ppbackr{i}{i}^*\cone{L}_{\group{S}}
  \end{tikzcd}
\]

\begin{proposition}\label{PropositionNuUnit}
  The action $\nu$ is unital.
\end{proposition}
\begin{proof}
  First, note that there is a diagonal morphism $\Delta \colon S \to
  \ppback{i}{i}$ satisfying $\ppbackl{i}{i}\Delta = \ppbackr{i}{i}\Delta = \id_{S}$
  and $\lambda\Delta = \id_{i}$. We thus have the following relation:
  \begin{equation}
    \begin{split}
     &\begin{tikzpicture}[baseline=(current bounding box.center)]
       \node[tip] (end1) at (0, -1) {$i_!$};
       \node[tip] (end2) at (1, -1) {$i^*$};
       \node[dot] (uni1) at (0.5, 0) {};
       \draw (uni1.center) to[out=180, in=90] (end1);
       \draw (uni1.center) to[out=0 , in=90] (end2);
     \end{tikzpicture} \spacedequal{}
     \begin{tikzpicture}[baseline=(current bounding box.center)]
       \node[tip] (end1) at (0, -3) {$i_!$};
       \node[tip] (end2) at (1, -3) {$i^*$};
       \node[dot] (uni1) at (0.5, 0) {};
       \node[iden] (id1) at (2.5, 0) {};
       \node[naturaltr] (lam) at (1.5, -1) {$\lambda^*$};
       \node[iden] (id2) at (2.5, -2) {};
       \draw (uni1.center) to[out=180, in=90] (end1);
       \draw (uni1.center) to[out=0 , in=90] (lam.north west);
       \draw (id1.south west) to[out=225, in=90] (lam.north east);
       \draw (id1.south east) to[out=315, in=45] node [right] {$\Delta^*$} (id2.north east);
       \draw (lam.south west) to[out=270, in=90] (end2);
       \draw (lam.south east) to[out=270, in=135] (id2.north west);
     \end{tikzpicture}  \\ &= \quad
     \begin{tikzpicture}[baseline=(current bounding box.center)]
       \node[tip] (end1) at (0, -3) {$i_!$};
       \node[tip] (end2) at (1, -3) {$i^*$};
       \node[dot] (uni1) at (0.5, 0) {};
       \node[iden] (id1) at (3.5, 0) {};
       \node[naturaltr] (lam) at (1.5, -1) {$\lambda^*$};
       \node[dot] (uni2) at (3, -1.5) {};
       \node[dot] (cou2) at (2.5, -2) {};
       \node[iden] (id2) at (3.5, -2) {};
       \draw (uni1.center) to[out=180, in=90] (end1);
       \draw (uni1.center) to[out=0 , in=90] (lam.north west);
       \draw (id1.south west) to[out=225, in=90] (lam.north east);
       \draw (id1.south east) to[out=315, in=45] node [right] {$\Delta^*$} (id2.north east);
       \draw (lam.south west) to[out=270, in=90] (end2);
       \draw (lam.south east) to[out=270, in=180] (cou2.center);
       \draw (uni2.center) to[out=180, in=0] (cou2.center);
       \draw (uni2.center) to[out=0, in=135] (id2.north west);
     \end{tikzpicture} \spacedequal{}
     \begin{tikzpicture}[baseline=(current bounding box.center)]
       \node[tip] (end1) at (0, -3) {$i_!$};
       \node[tip] (end2) at (1, -3) {$i^*$};
       \node[iden] (id1) at (2, 0) {};
       \node[dot] (uni2) at (1.5, -1) {};
       \node[iden] (id2) at (2, -2) {};
       \node[naturaltr] (bet) at (0.5, -2) {$\lambda_!$};
       \draw (id1.south west) to[out=225, in=90] (bet.north west);
       \draw (id1.south east) to[out=315, in=45] node [right] {$\Delta^*$} (id2.north east);
       \draw (uni2.center) to[out=180, in=90] (bet.north east);
       \draw (uni2.center) to[out=0, in=135] (id2.north west);
       \draw (bet.south west) to[out=270, in=90] (end1);
       \draw (bet.south east) to[out=270, in=90] (end2);
       \end{tikzpicture}
    \end{split}\label{eqUnitDelta}
  \end{equation}
   Hence, we have:
   
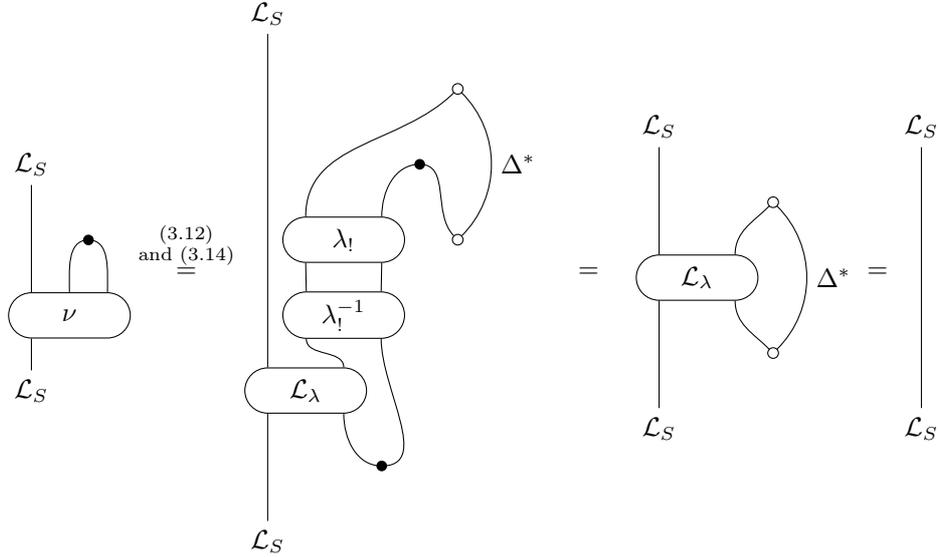
\begin{align*}
     \begin{tikzpicture}[baseline=(current bounding box.center)]
       \node[tip] (sta1) at (0, 0) {$\cone{L}_S$};
       \node[tip] (end1) at (0, -3) {$\cone{L}_S$};
       \node[dot] (uni1) at (0.75, -1) {};
       \node[naturaltr] (alp1) at (0.5, -2) {$\nu$};
       \draw (sta1) to[out=270, in=90] (alp1.north west);
       \draw (uni1.center) to[out=180, in=90] (alp1.north);
       \draw (uni1.center) to[out=0 , in=90] (alp1.north east);
       \draw (alp1.south west) to[out=270, in=90] (end1);
     \end{tikzpicture} \overset{\substack{\eqref{defNu} \\\text{and } \eqref{eqUnitDelta}}}{=}
     \begin{tikzpicture}[baseline=(current bounding box.center)]
       \node[tip] (sta1) at (0, 1) {$\cone{L}_S$};
       \node[tip] (end1) at (0, -6) {$\cone{L}_S$};
       \node[iden] (id1) at (2.5, 0) {};
       \node[dot] (uni2) at (2, -1) {};
       \node[iden] (id2) at (2.5, -2) {};
       \node[naturaltr] (bet1) at (1, -2) {$\lambda_!$};
       \node[naturaltr] (bet2) at (1, -3) {$\lambda_!^{-1}$};
       \node[naturaltr] (zet1) at (0.5, -4) {$\cone{L}_{\lambda}$};
       \node[dot] (cou1) at (1.5, -5) {};
       \draw (sta1) to[out=270, in=90] (zet1.north west);
       \draw (id1.south west) to[out=225, in=90] (bet1.north west);
       \draw (id1.south east) to[out=315, in=45] node [right] {$\Delta^*$} (id2.north east);
       \draw (uni2.center) to[out=180, in=90] (bet1.north east);
       \draw (uni2.center) to[out=0, in=135] (id2.north west);
       \draw (bet1.south west) to[out=270, in=90] (bet2.north west);
       \draw (bet1.south east) to[out=270, in=90] (bet2.north east);
       \draw (bet2.south west) to[out=270, in=90] (zet1.north east);
       \draw (bet2.south east) to[out=270, in=0]  (cou1.center);
       \draw (zet1.south west) to[out=270, in=90] (end1);
       \draw (zet1.south east) to[out=270, in=180] (cou1.center);
       \end{tikzpicture} \spacedequal{}
     \begin{tikzpicture}[baseline=(current bounding box.center)]
       \node[tip] (sta1) at (0, 1) {$\cone{L}_S$};
       \node[tip] (end1) at (0, -3) {$\cone{L}_S$};
       \node[iden] (id1) at (1.5, 0) {};
       \node[naturaltr] (zet1) at (0.5, -1) {$\cone{L}_{\lambda}$};
       \node[iden] (id2) at (1.5, -2) {};
       \draw (sta1) to[out=270, in=90] (zet1.north west);
       \draw (id1.south west) to[out=225, in=90] (zet1.north east);
       \draw (id1.south east) to[out=315, in=45] node [right] {$\Delta^*$} (id2.north east);
       \draw (zet1.south west) to[out=270, in=90] (end1);
       \draw (zet1.south east) to[out=270, in=135] (id2.north west);
     \end{tikzpicture} =
     \begin{tikzpicture}[baseline=(current bounding box.center)]
       \node[tip] (sta1) at (0, 1) {$\cone{L}_S$};
       \node[tip] (end1) at (0, -3) {$\cone{L}_S$};
       \draw (sta1) to[out=270, in=90] (end1);
     \end{tikzpicture}
   \end{align*}
\end{proof}

\begin{proposition}\label{PropositionNuAssociative}
  The action $\nu$ is associative.
\end{proposition}
\begin{proof}
  We introduce another 2-pullback square
  
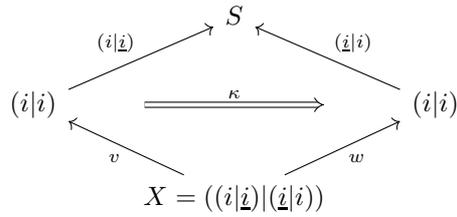
\[
    \begin{tikzcd}
      & S & \\
      \ppback{i}{i} \arrow[ru, "\ppbackr{i}{i}"]
      \arrow[rr, "\kappa", Rightarrow, shorten=10mm] & &
      \ppback{i}{i} \arrow[lu, "\ppbackl{i}{i}", swap] \\
      & X = \ppback{\ppbackr{i}{i}}{\ppbackl{i}{i}} \arrow[lu, "v"] \arrow[ru, "w", swap] &
    \end{tikzcd}
  \]
  The universal property of $\ppback{i}{i}$ allows us to define a morphism
  $\nabla \colon X \to \ppback{i}{i}$ such that
 \begin{equation}\label{EqNabla}
   \begin{tikzpicture}[baseline=(current bounding box.center)]
     \node[tip] (sta1) at (0,0)  {$v$};
     \node[tip] (sta2) at (1,0)  {$\ppbackl{i}{i}$};
     \node[tip] (sta3) at (2,0)  {$i$};
     \node[tip] (end1) at (0,-4) {$w$};
     \node[tip] (end2) at (1,-4) {$\ppbackr{i}{i}$};
     \node[tip] (end3) at (2,-4) {$i$};
     \node[naturaltr] (lam1) at (1.5 ,-1) {$\lambda$};
     \node[naturaltr] (kap)  at (0.5,-2) {$\kappa$};
     \node[naturaltr] (lam2) at (1.5 ,-3) {$\lambda$};
     \draw (sta2)     to[out=270, in=90] (lam1.north west);
     \draw (sta3)     to[out=270, in=90 ] (lam1.north east);
     \draw (sta1)     to[out=270, in=90] (kap.north west);
     \draw (lam1.south west) to[out=270, in=90 ] (kap.north east);
     \draw (kap.south east) to[out=270, in=90] (lam2.north west);
     \draw (lam1.south) to                  node[right] {$i$} (lam2.north);
     \draw (kap.south west) to[out=270, in=90] (end1);
     \draw (lam2.south west) to[out=270, in=90] (end2);
     \draw (lam2.south east) to[out=270, in=90] (end3);
   \end{tikzpicture}
   \spacedequal
   \begin{tikzpicture}[baseline=(current bounding box.center)]
     \node[tip] (sta1) at (0,0)  {$v$};
     \node[tip] (sta2) at (1,0)  {$\ppbackl{i}{i}$};
     \node[tip] (sta3) at (2,0)  {$i$};
     \node[tip] (end1) at (0,-4) {$w$};
     \node[tip] (end2) at (1,-4) {$\ppbackr{i}{i}$};
     \node[tip] (end3) at (2,-4) {$i$};
     \node[iden] (id1) at (0.75 ,-1) {};
     \node[naturaltr] (lam)  at (1.5  ,-2) {$\lambda$};
     \node[iden] (id2) at (0.75 ,-3) {};
     \draw (sta1)    to[out=270, in=135] (id1.north west);
     \draw (sta2)    to[out=270, in=45] (id1.north east);
     \draw (sta3)    to[out=270, in=90] (lam.north east);
     \draw (id1.south west) to[out=225, in=135] node[left] {$\nabla$} (id2.north west);
     \draw (id1.south east) to[out=315, in=90] (lam.north west);
     \draw (lam.south west) to[out=270, in=45] (id2.north east);
     \draw (id2.south west) to[out=225, in=90 ] (end1);
     \draw (id2.south east) to[out=315, in=90 ] (end2);
     \draw (lam.south east) to[out=270, in=90 ] (end3);
   \end{tikzpicture}
 \end{equation}

 Applying the 2-functor $\ffunctor{M}$ on \cref{EqNabla} yields
 \begin{equation}\label{EqNablaMack}
   \begin{tikzpicture}[baseline=(current bounding box.center)]
     \node[tip] (sta1) at (2,0)  {$v^*$};
     \node[tip] (sta2) at (1,0)  {$\ppbackl{i}{i}^*$};
     \node[tip] (sta3) at (0,0)  {$i^*$};
     \node[tip] (end1) at (2,-4) {$w^*$};
     \node[tip] (end2) at (1,-4) {$\ppbackr{i}{i}^*$};
     \node[tip] (end3) at (0,-4) {$i^*$};
     \node[naturaltr] (lam1) at (0.5 ,-1) {$\lambda^*$};
     \node[naturaltr] (kap)  at (1.5,-2) {$\kappa^*$};
     \node[naturaltr] (lam2) at (0.5 ,-3) {$\lambda^*$};
     \draw (sta2)     to[out=270, in=90 ] (lam1.north east);
     \draw (sta3)     to[out=270, in=90] (lam1.north west);
     \draw (sta1)     to[out=270, in=90 ] (kap.north east);
     \draw (lam1.south east) to[out=270, in=90] (kap.north west);
     \draw (kap.south west) to[out=270, in=90 ] (lam2.north east);
     \draw (lam1.south) to                  node[left]   {$i^*$} (lam2.north);
     \draw (kap.south east) to[out=270, in=90 ] (end1);
     \draw (lam2.south east) to[out=270, in=90 ] (end2);
     \draw (lam2.south west) to[out=270, in=90 ] (end3);
   \end{tikzpicture}
   \spacedequal
   \begin{tikzpicture}[baseline=(current bounding box.center)]
     \node[tip] (sta1) at (2,0)  {$v^*$};
     \node[tip] (sta2) at (1,0)  {$\ppbackl{i}{i}^*$};
     \node[tip] (sta3) at (0,0)  {$i^*$};
     \node[tip] (end1) at (2,-4) {$w^*$};
     \node[tip] (end2) at (1,-4) {$\ppbackr{i}{i}^*$};
     \node[tip] (end3) at (0,-4) {$i^*$};
     \node[iden] (id1) at (1.25 ,-1) {};
     \node[naturaltr] (lam)  at (0.5  ,-2) {$\lambda^*$};
     \node[iden] (id2) at (1.25 ,-3) {};
     \draw (sta1)    to[out=270, in=45 ] (id1.north east);
     \draw (sta2)    to[out=270, in=135] (id1.north west);
     \draw (sta3)    to[out=270, in=90] (lam.north west);
     \draw (id1.south east) to[out=315, in=45] node[right] {$\nabla^*$} (id2.north east);
     \draw (id1.south west) to[out=225, in=90] (lam.north east);
     \draw (lam.south east) to[out=270, in=135] (id2.north west);
     \draw (id2.south east) to[out=315, in=90 ] (end1);
     \draw (id2.south west) to[out=225, in=90 ] (end2);
     \draw (lam.south west) to[out=270, in=90 ] (end3);
   \end{tikzpicture}
 \end{equation}
Similarly, the compatibility between the cone \cone{L} and the 2-morphisms of
\cattransportext{\group{G}}{\group{S}} provides the equation
 \begin{equation}\label{eqNabla2Lim}
   \begin{tikzpicture}[baseline=(current bounding box.center)]
     \node[tip] (sta1) at (2,0)  {$v^*$};
     \node[tip] (sta2) at (1,0)  {$\ppbackl{i}{i}^*$};
     \node[tip] (sta3) at (0,0)  {$\cone{L}_{\group{S}}$};
     \node[tip] (end1) at (2,-4) {$w^*$};
     \node[tip] (end2) at (1,-4) {$\ppbackr{i}{i}^*$};
     \node[tip] (end3) at (0,-4) {$\cone{L}_{\group{S}}$};
     \node[naturaltr] (lam1) at (0.5 ,-1) {$\cone{L}_{\lambda}$};
     \node[naturaltr] (kap)  at (1.5,-2) {$\kappa^*$};
     \node[naturaltr] (lam2) at (0.5 ,-3) {$\cone{L}_{\lambda}$};
     \draw (sta2)     to[out=270, in=90 ] (lam1.north east);
     \draw (sta3)     to[out=270, in=90] (lam1.north west);
     \draw (sta1)     to[out=270, in=90 ] (kap.north east);
     \draw (lam1.south east) to[out=270, in=90] (kap.north west);
     \draw (kap.south west) to[out=270, in=90 ] (lam2.north east);
     \draw (lam1.south) to                  node[left]   {$\cone{L}_{\group{S}}$} (lam2.north);
     \draw (kap.south east) to[out=270, in=90 ] (end1);
     \draw (lam2.south east) to[out=270, in=90 ] (end2);
     \draw (lam2.south west) to[out=270, in=90 ] (end3);
   \end{tikzpicture}
   \spacedequal
   \begin{tikzpicture}[baseline=(current bounding box.center)]
     \node[tip] (sta1) at (2,0)  {$v^*$};
     \node[tip] (sta2) at (1,0)  {$\ppbackl{i}{i}^*$};
     \node[tip] (sta3) at (0,0)  {$\cone{L}_{\group{S}}$};
     \node[tip] (end1) at (2,-4) {$w^*$};
     \node[tip] (end2) at (1,-4) {$\ppbackr{i}{i}^*$};
     \node[tip] (end3) at (0,-4) {$\cone{L}_{\group{S}}$};
     \node[iden] (id1) at (1.25 ,-1) {};
     \node[naturaltr] (lam)  at (0.5  ,-2) {$\cone{L}_{\lambda}$};
     \node[iden] (id2) at (1.25 ,-3) {};
     \draw (sta1)    to[out=270, in=45 ] (id1.north east);
     \draw (sta2)    to[out=270, in=135] (id1.north west);
     \draw (sta3)    to[out=270, in=90] (lam.north west);
     \draw (id1.south east) to[out=315, in=45] node[right] {$\nabla^*$} (id2.north east);
     \draw (id1.south west) to[out=225, in=90] (lam.north east);
     \draw (lam.south east) to[out=270, in=135] (id2.north west);
     \draw (id2.south east) to[out=315, in=90 ] (end1);
     \draw (id2.south west) to[out=225, in=90 ] (end2);
     \draw (lam.south west) to[out=270, in=90 ] (end3);
   \end{tikzpicture}
 \end{equation}
 We can now prove the associativity of $\nu$. We want to prove :
   \begin{equation}
     \label{EqNuAssoc}
     \begin{tikzpicture}[baseline=(current bounding box.center)]
       \node[tip] (sta1) at (0, 0) {$\cone{L}_S$};
       \node[tip] (sta2) at (0.5, 0) {$i_!$};
       \node[tip] (sta3) at (1, 0) {$i^*$};
       \node[tip] (sta4) at (1.5, 0) {$i_!$};
       \node[tip] (sta5) at (2, 0) {$i^*$};
       \node[tip] (end1) at (1.5, -3) {$\cone{L}_S$};
       \node[naturaltr] (alp1) at (0.5, -1) {$\nu$};
       \node[naturaltr] (alp2) at (1.5, -2) {$\nu$};
       \draw (sta1) to[out=270, in=90] (alp1.north west);
       \draw (sta2) to[out=270, in=90] (alp1.north);
       \draw (sta3) to[out=270, in=90] (alp1.north east);
       \draw (sta4) to[out=270, in=90] (alp2.north);
       \draw (sta5) to[out=270, in=90] (alp2.north east);
       \draw (alp1.south east) to[out=270, in=90] (alp2.north west);
       \draw (alp2.south) to[out=270, in=90] (end1);
     \end{tikzpicture} \spacedequal{}
      \begin{tikzpicture}[baseline=(current bounding box.center)]
       \node[tip] (sta1) at (0, 0) {$\cone{L}_S$};
       \node[tip] (sta2) at (0.5, 0) {$i_!$};
       \node[tip] (sta3) at (1, 0) {$i^*$};
       \node[tip] (sta4) at (1.5, 0) {$i_!$};
       \node[tip] (sta5) at (2, 0) {$i^*$};
       \node[tip] (end1) at (0.5, -3) {$\cone{L}_S$};
       \node[dot] (cou1) at (1.25, -1) {};
       \node[naturaltr] (alp2) at (0.5, -2) {$\nu$};
       \draw (sta1) to[out=270, in=90] (alp2.north west);
       \draw (sta2) to[out=270, in=90] (alp2.north);
       \draw (sta3) to[out=270, in=180] (cou1.center);
       \draw (sta4) to[out=270, in=0] (cou1.center);
       \draw (sta5) to[out=270, in=90] (alp2.north east);
       \draw (alp2.south) to[out=270, in=90] (end1);
     \end{tikzpicture}
   \end{equation}
   Consider the mate $\kappa_!$ of
   $\kappa^*$:
   \[
     \kappa_!=
      \begin{tikzpicture}[baseline=(current bounding box.center)]
       \node[tip] (sta1) at (2, 0) {$v^*$};
       \node[tip] (sta2) at (3, 0) {$w_!$};
       \node[tip] (end1) at (0, -4) {$\ppbackr{i}{i}_!$};
       \node[tip] (end2) at (1, -4) {$\ppbackl{i}{i}^*$};
       \node[dot] (uni1) at (0.5, -1) {};
       \node[naturaltr] (two) at (1.5, -2) {$\kappa^*$};
       \node[dot] (cou1) at (2.5, -3) {};
       \draw (sta1) to[out=270, in=90] (two.north east);
       \draw (sta2) to[out=270, in=0] (cou1.center);
       \draw (uni1.center) to[out=180, in=90] (end1);
       \draw (uni1.center) to[out=0, in=90] (two.north west);
       \draw (two.south west) to[out=270, in=90] (end2);
       \draw (two.south east) to[out=270, in=180] (cou1.center);
     \end{tikzpicture}
   \]
   Since $\lambda_!$ and $\kappa_!$ are invertible, \cref{EqNuAssoc} is equivalent to the equation:
   \begin{equation}\label{EqNuAssoc2}
     \begin{tikzpicture}[baseline=(current bounding box.center)]
       \node[tip] (sta1) at (0.5, 0) {$\cone{L}_S$};
       \node[tip] (sta2) at (1.25, 0) {$\ppbackl{i}{i}^*$};
       \node[tip] (sta3) at (2, 0) {$v^*$};
       \node[tip] (sta4) at (3, 0) {$w_!$};
       \node[tip] (sta5) at (4, 0) {$\ppbackr{i}{i}_!$};
       \node[tip] (end1) at (3, -5) {$\cone{L}_S$};
       \node[naturaltr] (gam1) at (2.5, -1) {$\kappa_!$};
       \node[naturaltr] (bet1) at (1.5, -2) {$\lambda_!$};
       \node[naturaltr] (bet2) at (3.5, -2) {$\lambda_!$};
       \node[naturaltr] (alp1) at (1, -3) {$\nu$};
       \node[naturaltr] (alp2) at (3, -4) {$\nu$};
       \draw (sta1) to[out=270, in=90] (alp1.north west);
       \draw (sta2) to[out=270, in=90] (bet1.north west);
       \draw (sta3) to[out=270, in=90] (gam1.north west);
       \draw (sta4) to[out=270, in=90] (gam1.north east);
       \draw (sta5) to[out=270, in=90] (bet2.north east);
       \draw (gam1.south west) to[out=270, in=90] (bet1.north east);
       \draw (gam1.south east) to[out=270, in=90] (bet2.north west);
       \draw (bet1.south west) to[out=270, in=90] (alp1.north);
       \draw (bet1.south) to[out=270, in=90] (alp1.north east);
       \draw (bet2.south west) to[out=270, in=90] (alp2.north);
       \draw (bet2.south) to[out=270, in=90] (alp2.north east);
       \draw (alp1.south east) to[out=270, in=90] (alp2.north west);
       \draw (alp2.south) to[out=270, in=90] (end1);
     \end{tikzpicture} \spacedequal{}
      \begin{tikzpicture}[baseline=(current bounding box.center)]
       \node[tip] (sta1) at (0.5, 0) {$\cone{L}_S$};
       \node[tip] (sta2) at (1.25, 0) {$\ppbackl{i}{i}^*$};
       \node[tip] (sta3) at (2, 0) {$v^*$};
       \node[tip] (sta4) at (3, 0) {$w_!$};
       \node[tip] (sta5) at (4, 0) {$\ppbackr{i}{i}_!$};
       \node[tip] (end1) at (1, -5) {$\cone{L}_S$};
       \node[naturaltr] (gam1) at (2.5, -1) {$\kappa_!$};
       \node[naturaltr] (bet1) at (1.5, -2) {$\lambda_!$};
       \node[naturaltr] (bet2) at (3.5, -2) {$\lambda_!$};
       \node[dot] (cou1) at (2.25, -2.75) {};
       \node[naturaltr] (alp2) at (1, -4) {$\nu$};
       \draw (sta1) to[out=270, in=90] (alp2.north west);
       \draw (sta2) to[out=270, in=90] (bet1.north west);
       \draw (sta3) to[out=270, in=90] (gam1.north west);
       \draw (sta4) to[out=270, in=90] (gam1.north east);
       \draw (sta5) to[out=270, in=90] (bet2.north east);
       \draw (gam1.south west) to[out=270, in=90] (bet1.north east);
       \draw (gam1.south east) to[out=270, in=90] (bet2.north west);
       \draw (bet1.south west) to[out=270, in=90] (alp2.north);
       \draw (bet1.south) to[out=270, in=180] (cou1.center);
       \draw (bet2.south west) to[out=270, in=0] (cou1.center);
       \draw (bet2.south) to[out=270, in=90] (alp2.north east);
       \draw (alp2.south) to[out=270, in=90] (end1);
     \end{tikzpicture}
   \end{equation}
   Expanding the definitions of $\nu$ and $\kappa_!$ on the left-hand
   side of \cref{EqNuAssoc2} leads to the following computation:
   \begin{align*}
     &\begin{tikzpicture}[baseline=(current bounding box.center)]
       \node[tip] (sta1) at (0.5, 0) {$\cone{L}_S$};
       \node[tip] (sta2) at (1.25, 0) {$\ppbackl{i}{i}^*$};
       \node[tip] (sta3) at (2, 0) {$v^*$};
       \node[tip] (sta4) at (3, 0) {$w_!$};
       \node[tip] (sta5) at (4, 0) {$\ppbackr{i}{i}_!$};
       \node[tip] (end1) at (3, -5) {$\cone{L}_S$};
       \node[naturaltr] (gam1) at (2.5, -1) {$\kappa_!$};
       \node[naturaltr] (bet1) at (1.5, -2) {$\lambda_!$};
       \node[naturaltr] (bet2) at (3.5, -2) {$\lambda_!$};
       \node[naturaltr] (alp1) at (1, -3) {$\nu$};
       \node[naturaltr] (alp2) at (3, -4) {$\nu$};
       \draw (sta1) to[out=270, in=90] (alp1.north west);
       \draw (sta2) to[out=270, in=90] (bet1.north west);
       \draw (sta3) to[out=270, in=90] (gam1.north west);
       \draw (sta4) to[out=270, in=90] (gam1.north east);
       \draw (sta5) to[out=270, in=90] (bet2.north east);
       \draw (gam1.south west) to[out=270, in=90] (bet1.north east);
       \draw (gam1.south east) to[out=270, in=90] (bet2.north west);
       \draw (bet1.south west) to[out=270, in=90] (alp1.north);
       \draw (bet1.south) to[out=270, in=90] (alp1.north east);
       \draw (bet2.south west) to[out=270, in=90] (alp2.north);
       \draw (bet2.south) to[out=270, in=90] (alp2.north east);
       \draw (alp1.south east) to[out=270, in=90] (alp2.north west);
       \draw (alp2.south) to[out=270, in=90] (end1);
     \end{tikzpicture} \xspacedequal{\eqref{defNu}}
      \begin{tikzpicture}[baseline=(current bounding box.center)]
       \node[tip] (sta1) at (0, 0) {$\cone{L}_S$};
       \node[tip] (sta2) at (1, 0) {$\ppbackl{i}{i}^*$};
       \node[tip] (sta3) at (2, 0) {$v^*$};
       \node[tip] (sta4) at (3, 0) {$w_!$};
       \node[tip] (sta5) at (4, 0) {$\ppbackr{i}{i}_!$};
       \node[tip] (end1) at (2, -8) {$\cone{L}_S$};
       \node[naturaltr] (gam1) at (2.5, -1) {$\kappa_!$};
       \node[naturaltr] (bet1) at (1.5, -2) {$\lambda_!$};
       \node[naturaltr] (bet2) at (3.5, -2) {$\lambda_!$};
       \node[naturaltr] (bet3) at (1.5, -3) {$\lambda_!^{-1}$};
       \node[naturaltr] (bet4) at (3.5, -3) {$\lambda_!^{-1}$};
       \node[naturaltr] (zet1) at (0.5, -4) {$\cone{L}_{\lambda}$};
       \node[dot] (cou1) at (1.5, -4.75) {};
       \node[naturaltr] (zet2) at (2.5, -6) {$\cone{L}_{\lambda}$};
       \node[dot] (cou2) at (3.5, -6.75) {};
       \draw (sta1) to[out=270, in=90] (zet1.north west);
       \draw (sta2) to[out=270, in=90] (bet1.north west);
       \draw (sta3) to[out=270, in=90] (gam1.north west);
       \draw (sta4) to[out=270, in=90] (gam1.north east);
       \draw (sta5) to[out=270, in=90] (bet2.north east);
       \draw (gam1.south west) to[out=270, in=90] (bet1.north east);
       \draw (gam1.south east) to[out=270, in=90] (bet2.north west);
       \draw (bet1.south west) to[out=270, in=90] (bet3.north west);
       \draw (bet1.south east) to[out=270, in=90] (bet3.north east);
       \draw (bet2.south west) to[out=270, in=90] (bet4.north west);
       \draw (bet2.south east) to[out=270, in=90] (bet4.north east);
       \draw (bet3.south west) to[out=270, in=90] (zet1.north east);
       \draw (bet3.south east) to[out=270, in=0] (cou1.center);
       \draw (bet4.south west) to[out=270, in=90] (zet2.north east);
       \draw (bet4.south east) to[out=270, in=0] (cou2.center);
       \draw (zet1.south west) to[out=270, in=90] (zet2.north west);
       \draw (zet1.south east) to[out=270, in=180] (cou1.center);
       \draw (zet2.south west) to[out=270, in=90] (end1);
       \draw (zet2.south east) to[out=270, in=180] (cou2.center);
     \end{tikzpicture} \\ &\overset{\eqref{eqDefMate}}{=} \quad
      \begin{tikzpicture}[baseline=(current bounding box.center)]
       \node[tip] (sta1) at (0, 0) {$\cone{L}_S$};
       \node[tip] (sta2) at (1, 0) {$\ppbackl{i}{i}^*$};
       \node[tip] (sta3) at (4, 0) {$v^*$};
       \node[tip] (sta4) at (5, 0) {$w_!$};
       \node[tip] (sta5) at (6, 0) {$\ppbackr{i}{i}_!$};
       \node[tip] (end1) at (2, -4) {$\cone{L}_S$};
       \node[dot] (uni1) at (2.5, -0.5) {};
       \node[naturaltr] (two) at (3.5, -1) {$\kappa^*$};
       \node[dot] (cou3) at (4.5, -1.5) {};
       \node[naturaltr] (zet1) at (0.5, -1) {$\cone{L}_{\lambda}$};
       \node[dot] (cou1) at (1.5, -1.5) {};
       \node[naturaltr] (zet2) at (2.5, -3) {$\cone{L}_{\lambda}$};
       \node[dot] (cou2) at (3.5, -3.5) {};
       \draw (sta1) to[out=270, in=90] (zet1.north west);
       \draw (sta2) to[out=270, in=90] (zet1.north east);
       \draw (sta3) to[out=270, in=90] (two.north east);
       \draw (sta4) to[out=270, in=0] (cou3.center);
       \draw (sta5) to[out=270, in=0] (cou2.center);
       \draw (uni1.center) to[out=180, in=0] (cou1.center);
       \draw (uni1.center) to[out=0, in=90] (two.north west);
       \draw (two.south west) to[out=270, in=90] (zet2.north east);
       \draw (two.south east) to[out=270, in=180] (cou3.center);
       \draw (zet1.south west) to[out=270, in=90] (zet2.north west);
       \draw (zet1.south east) to[out=270, in=180] (cou1.center);
       \draw (zet2.south west) to[out=270, in=90] (end1);
       \draw (zet2.south east) to[out=270, in=180] (cou2.center);
     \end{tikzpicture}   \\
     & = \quad
      \begin{tikzpicture}[baseline=(current bounding box.center)]
       \node[tip] (sta1) at (0, 0) {$\cone{L}_S$};
       \node[tip] (sta2) at (1, 0) {$\ppbackl{i}{i}^*$};
       \node[tip] (sta3) at (2, 0) {$v^*$};
       \node[tip] (sta4) at (3, 0) {$w_!$};
       \node[tip] (sta5) at (4, 0) {$\ppbackr{i}{i}_!$};
       \node[tip] (end1) at (0, -5) {$\cone{L}_S$};
       \node[naturaltr] (zet1) at (0.5, -1) {$\cone{L}_{\lambda}$};
       \node[naturaltr] (two) at (1.5, -2) {$\kappa^*$};
       \node[dot] (cou3) at (2.5, -3) {};
       \node[naturaltr] (zet2) at (0.5, -3) {$\cone{L}_{\lambda}$};
       \node[dot] (cou2) at (1.5, -4) {};
       \draw (sta1) to[out=270, in=90] (zet1.north west);
       \draw (sta2) to[out=270, in=90] (zet1.north east);
       \draw (sta3) to[out=270, in=90] (two.north east);
       \draw (sta4) to[out=270, in=0] (cou3.center);
       \draw (sta5) to[out=270, in=0] (cou2.center);
       \draw (zet1.south west) to[out=270, in=90] (zet2.north west);
       \draw (zet1.south east) to[out=270, in=90] (two.north west);
       \draw (two.south west) to[out=270, in=90] (zet2.north east);
       \draw (two.south east) to[out=270, in=180] (cou3.center);
       \draw (zet2.south west) to[out=270, in=90] (end1);
       \draw (zet2.south east) to[out=270, in=180] (cou2.center);
     \end{tikzpicture} \overset{\eqref{eqNabla2Lim}}{=} \quad
       \begin{tikzpicture}[baseline=(current bounding box.center)]
       \node[tip] (sta1) at (0, 0) {$\cone{L}_S$};
       \node[tip] (sta2) at (1, 0) {$\ppbackl{i}{i}^*$};
       \node[tip] (sta3) at (2, 0) {$v^*$};
       \node[tip] (sta4) at (3, 0) {$w_!$};
       \node[tip] (sta5) at (4, 0) {$\ppbackr{i}{i}_!$};
       \node[tip] (end1) at (0, -5) {$\cone{L}_S$};
       \node[iden] (one) at (1.5, -1) {};
       \node[naturaltr] (zet) at (0.5, -2) {$\cone{L}_{\lambda}$};
       \node[iden] (thr) at (1.5, -3) {};
       \node[dot] (cou3) at (2.5, -3.5) {};
       \node[dot] (cou2) at (1.5, -4.5) {};
       \draw (sta1) to[out=270, in=90] (zet.north west);
       \draw (sta2) to[out=270, in=90] (one.north west);
       \draw (sta3) to[out=270, in=90] (one.north east);
       \draw (sta4) to[out=270, in=0] (cou3.center);
       \draw (sta5) to[out=270, in=0] (cou2.center);
       \draw (one.south west) to[out=270, in=90] (zet.north east);
       \draw (one.south east) to[out=315, in=45] node[right]{$\nabla^*$} (thr.north east);
       \draw (zet.south west) to[out=270, in=90] (end1);
       \draw (zet.south east) to[out=270, in=90] (thr.north west);
       \draw (thr.south west) to[out=270, in=180] (cou2.center);
       \draw (thr.south east) to[out=270, in=180] (cou3.center);
     \end{tikzpicture}
   \end{align*}
   To simplify the right-hand side of \cref{EqNuAssoc2}, first note that we have the following relation:
   \begin{equation} \label{eqAlphaZetaHat}
     \begin{split}
       \begin{tikzpicture}[baseline=(current bounding box.center)]
         \node[tip] (sta1) at (0, 0) {$\cone{L}_S$};
         \node[tip] (sta2) at (2, 0) {$\ppbackl{i}{i}^*$};
         \node[tip] (end1) at (0, -3) {$\cone{L}_S$};
         \node[tip] (end2) at (2, -3) {$\ppbackr{i}{i}^*$};
         \node[dot] (uni) at (0.75, -0.5) {};
         \node[naturaltr] (lam) at (1.5, -1) {$\lambda^*$};
         \node[naturaltr] (alp) at (0.5, -2) {$\nu$};
         \draw (sta1) to[out=270, in=90] (alp.north west);
         \draw (sta2) to[out=270, in=90] (lam.north east);
         \draw (uni.center) to[out=180, in=90] (alp.north);
         \draw (uni.center) to[out=0, in=90] (lam.north west);
         \draw (lam.south west) to[out=270, in=90] (alp.north east);
         \draw (lam.south east) to[out=270, in=90] (end2);
         \draw (alp.south west) to[out=270, in=90] (end1);
       \end{tikzpicture} &\spacedequal{}
       \begin{tikzpicture}[baseline=(current bounding box.center)]
         \node[tip] (sta1) at (0, 0) {$\cone{L}_S$};
         \node[tip] (sta2) at (0.75, 0) {$\ppbackl{i}{i}^*$};
         \node[tip] (end1) at (0, -3) {$\cone{L}_S$};
         \node[tip] (end2) at (2.5, -3) {$\ppbackr{i}{i}^*$};
         \node[dot] (uni) at (2, -0.5) {};
         \node[naturaltr] (bet) at (1, -1) {$\lambda_!$};
         \node[naturaltr] (alp) at (0.5, -2) {$\nu$};
         \draw (sta1) to[out=270, in=90] (alp.north west);
         \draw (sta2) to[out=270, in=90] (bet.north west);
         \draw (uni.center) to[out=180, in=90] (bet.north east);
         \draw (uni.center) to[out=0, in=90] (end2);
         \draw (bet.south west) to[out=270, in=90] (alp.north);
         \draw (bet.south) to[out=270, in=90] (alp.north east);
         \draw (alp.south west) to[out=270, in=90] (end1);
       \end{tikzpicture} \\ &\xspacedequal{\eqref{defNu}}
       \begin{tikzpicture}[baseline=(current bounding box.center)]
         \node[tip] (sta1) at (0, 0) {$\cone{L}_S$};
         \node[tip] (sta2) at (1, 0) {$\ppbackl{i}{i}^*$};
         \node[tip] (end1) at (0, -3) {$\cone{L}_S$};
         \node[tip] (end2) at (2.5, -3) {$\ppbackr{i}{i}^*$};
         \node[dot] (uni) at (2, -0.5) {};
         \node[naturaltr] (zet) at (0.5, -1.5) {$\cone{L}_{\lambda}$};
         \node[dot] (cou) at (1.5, -2.5) {};
         \draw (sta1) to[out=270, in=90] (zet.north west);
         \draw (sta2) to[out=270, in=90] (zet.north east);
         \draw (uni.center) to[out=180, in=0] (cou.center);
         \draw (uni.center) to[out=0, in=90] (end2);
         \draw (zet.south west) to[out=270, in=90] (end1);
         \draw (zet.south east) to[out=270, in=180] (cou.center);
       \end{tikzpicture} \spacedequal{}
       \begin{tikzpicture}[baseline=(current bounding box.center)]
         \node[tip] (sta1) at (0, 0) {$\cone{L}_S$};
         \node[tip] (sta2) at (1, 0) {$\ppbackl{i}{i}^*$};
         \node[tip] (end1) at (0, -3) {$\cone{L}_S$};
         \node[tip] (end2) at (1, -3) {$\ppbackr{i}{i}^*$};
         \node[naturaltr] (zet) at (0.5, -1.5) {$\cone{L}_{\lambda}$};
         \draw (sta1) to[out=270, in=90] (zet.north west);
         \draw (sta2) to[out=270, in=90] (zet.north east);
         \draw (zet.south west) to[out=270, in=90] (end1);
         \draw (zet.south east) to[out=270, in=90] (end2);
       \end{tikzpicture}
     \end{split}
   \end{equation}
   Then, by expanding the definitions of $\lambda_!$ and $\kappa_!$ in the right-hand
   side of \cref{EqNuAssoc2}, we get:
   \begin{align*}
     &\begin{tikzpicture}[baseline=(current bounding box.center)]
       \node[tip] (sta1) at (0.5, 0) {$\cone{L}_S$};
       \node[tip] (sta2) at (1.25, 0) {$\ppbackl{i}{i}^*$};
       \node[tip] (sta3) at (2, 0) {$v^*$};
       \node[tip] (sta4) at (3, 0) {$w_!$};
       \node[tip] (sta5) at (4, 0) {$\ppbackr{i}{i}_!$};
       \node[tip] (end1) at (1, -5) {$\cone{L}_S$};
       \node[naturaltr] (gam1) at (2.5, -1) {$\kappa_!$};
       \node[naturaltr] (bet1) at (1.5, -2) {$\lambda_!$};
       \node[naturaltr] (bet2) at (3.5, -2) {$\lambda_!$};
       \node[dot] (cou1) at (2.25, -2.75) {};
       \node[naturaltr] (alp2) at (1, -4) {$\nu$};
       \draw (sta1) to[out=270, in=90] (alp2.north west);
       \draw (sta2) to[out=270, in=90] (bet1.north west);
       \draw (sta3) to[out=270, in=90] (gam1.north west);
       \draw (sta4) to[out=270, in=90] (gam1.north east);
       \draw (sta5) to[out=270, in=90] (bet2.north east);
       \draw (gam1.south west) to[out=270, in=90] (bet1.north east);
       \draw (gam1.south east) to[out=270, in=90] (bet2.north west);
       \draw (bet1.south west) to[out=270, in=90] (alp2.north);
       \draw (bet1.south) to[out=270, in=180] (cou1.center);
       \draw (bet2.south west) to[out=270, in=0] (cou1.center);
       \draw (bet2.south) to[out=270, in=90] (alp2.north east);
       \draw (alp2.south) to[out=270, in=90] (end1);
     \end{tikzpicture} \overset{\eqref{eqDefMate}}{=}
     \begin{tikzpicture}[baseline=(current bounding box.center)]
       \node[tip] (sta1) at (0.5, 0) {$\cone{L}_S$};
       \node[tip] (sta2) at (3, 0) {$\ppbackl{i}{i}^*$};
       \node[tip] (sta3) at (5.5, 0) {$v^*$};
       \node[tip] (sta4) at (6.5, 0) {$w_!$};
       \node[tip] (sta5) at (7.5, 0) {$\ppbackr{i}{i}_!$};
       \node[tip] (end1) at (1, -5) {$\cone{L}_S$};
       \node[dot] (uni1) at (4, -0.5) {};
       \node[naturaltr] (two) at (5, -1) {$\kappa^*$};
       \node[dot] (cou1) at (6, -1.5) {};
       \node[dot] (uni2) at (1.5, -0.5) {};
       \node[naturaltr] (lam1) at (2.5, -1) {$\lambda^*$};
       \node[dot] (cou2) at (3.5, -1.5) {};
       \node[dot] (uni3) at (3, -2) {};
       \node[naturaltr] (lam2) at (4, -2.5) {$\lambda^*$};
       \node[dot] (cou3) at (5, -3) {};
       \node[dot] (cou4) at (2.5, -3) {};
       \node[naturaltr] (alp2) at (1, -4) {$\nu$};
       \draw (sta1) to[out=270, in=90] (alp2.north west);
       \draw (sta2) to[out=270, in=90] (lam1.north east);
       \draw (sta3) to[out=270, in=90] (two.north east);
       \draw (sta4) to[out=270, in=0] (cou1.center);
       \draw (sta5) to[out=270, in=0] (cou3.center);
       \draw (uni1.center) to[out=180, in=0] (cou2.center);
       \draw (uni1.center) to[out=0, in=90] (two.north west);
       \draw (two.south west) to[out=270, in=90] (lam2.north east);
       \draw (two.south east) to[out=270, in=180] (cou1.center);
       \draw (uni2.center) to[out=180, in=90] (alp2.north);
       \draw (uni2.center) to[out=0, in=90] (lam1.north west);
       \draw (lam1.south west) to[out=270, in=180] (cou4.center);
       \draw (lam1.south east) to[out=270, in=180] (cou2.center);
       \draw (uni3.center) to[out=180, in=0] (cou4.center);
       \draw (uni3.center) to[out=0, in=90] (lam2.north west);
       \draw (lam2.south west) to[out=270, in=90] (alp2.north east);
       \draw (lam2.south east) to[out=270, in=180] (cou3.center);
       \draw (alp2.south) to[out=270, in=90] (end1);
     \end{tikzpicture} \\ &\spacedequal{}
     \begin{tikzpicture}[baseline=(current bounding box.center)]
       \node[tip] (sta1) at (0.5, 0) {$\cone{L}_S$};
       \node[tip] (sta2) at (3, 0) {$\ppbackl{i}{i}^*$};
       \node[tip] (sta3) at (4, 0) {$v^*$};
       \node[tip] (sta4) at (5, 0) {$w_!$};
       \node[tip] (sta5) at (6, 0) {$\ppbackr{i}{i}_!$};
       \node[tip] (end1) at (1, -5) {$\cone{L}_S$};
       \node[dot] (uni2) at (1.5, -0.5) {};
       \node[naturaltr] (lam1) at (2.5, -1) {$\lambda^*$};
       \node[naturaltr] (two) at (3.5, -2) {$\kappa^*$};
       \node[dot] (cou1) at (4.5, -2.5) {};
       \node[naturaltr] (lam2) at (2.5, -3) {$\lambda^*$};
       \node[dot] (cou3) at (3.5, -3.5) {};
       \node[naturaltr] (alp2) at (1, -4) {$\nu$};
       \draw (sta1) to[out=270, in=90] (alp2.north west);
       \draw (sta2) to[out=270, in=90] (lam1.north east);
       \draw (sta3) to[out=270, in=90] (two.north east);
       \draw (sta4) to[out=270, in=0] (cou1.center);
       \draw (sta5) to[out=270, in=0] (cou3.center);
       \draw (uni2.center) to[out=180, in=90] (alp2.north);
       \draw (uni2.center) to[out=0, in=90] (lam1.north west);
       \draw (lam1.south west) to[out=270, in=90] (lam2.north west);
       \draw (lam1.south east) to[out=270, in=90] (two.north west);
       \draw (two.south west) to[out=270, in=90] (lam2.north east);
       \draw (two.south east) to[out=270, in=180] (cou1.center);
       \draw (lam2.south west) to[out=270, in=90] (alp2.north east);
       \draw (lam2.south east) to[out=270, in=180] (cou3.center);
       \draw (alp2.south) to[out=270, in=90] (end1);
     \end{tikzpicture} \\ &\overset{\eqref{EqNablaMack}}{=}
     \begin{tikzpicture}[baseline=(current bounding box.center)]
       \node[tip] (sta1) at (0.5, 0) {$\cone{L}_S$};
       \node[tip] (sta2) at (3, 0) {$\ppbackl{i}{i}^*$};
       \node[tip] (sta3) at (4, 0) {$v^*$};
       \node[tip] (sta4) at (5, 0) {$w_!$};
       \node[tip] (sta5) at (6, 0) {$\ppbackr{i}{i}_!$};
       \node[tip] (end1) at (1, -5) {$\cone{L}_S$};
       \node[dot] (uni2) at (1.5, -1.5) {};
       \node[iden] (one) at (3.5, -1) {};
       \node[naturaltr] (lam) at (2.5, -2) {$\lambda^*$};
       \node[iden] (thr) at (3.5, -3) {};
       \node[dot] (cou1) at (4.5, -3.5) {};
       \node[dot] (cou3) at (3.5, -4.5) {};
       \node[naturaltr] (alp2) at (1, -4) {$\nu$};
       \draw (sta1) to[out=270, in=90] (alp2.north west);
       \draw (sta2) to[out=270, in=90] (one.north west);
       \draw (sta3) to[out=270, in=90] (one.north east);
       \draw (sta4) to[out=270, in=0] (cou1.center);
       \draw (sta5) to[out=270, in=0] (cou3.center);
       \draw (uni2.center) to[out=180, in=90] (alp2.north);
       \draw (uni2.center) to[out=0, in=90] (lam.north west);
       \draw (one.south west) to[out=270, in=90] (lam.north east);
       \draw (one.south east) to[out=315, in=45] node[right]{$\nabla^*$} (thr.north east);
       \draw (lam.south west) to[out=270, in=90] (alp2.north east);
       \draw (lam.south east) to[out=270, in=90] (thr.north west);
       \draw (thr.south west) to[out=270, in=180] (cou3.center);
       \draw (thr.south east) to[out=270, in=180] (cou1.center);
       \draw (alp2.south) to[out=270, in=90] (end1);
     \end{tikzpicture} \overset{\eqref{eqAlphaZetaHat}}{=}
       \begin{tikzpicture}[baseline=(current bounding box.center)]
       \node[tip] (sta1) at (0, 0) {$\cone{L}_S$};
       \node[tip] (sta2) at (1, 0) {$\ppbackl{i}{i}^*$};
       \node[tip] (sta3) at (2, 0) {$v^*$};
       \node[tip] (sta4) at (3, 0) {$w_!$};
       \node[tip] (sta5) at (4, 0) {$\ppbackr{i}{i}_!$};
       \node[tip] (end1) at (0, -5) {$\cone{L}_S$};
       \node[iden] (one) at (1.5, -1) {};
       \node[naturaltr] (zet) at (0.5, -2) {$\cone{L}_{\lambda}$};
       \node[iden] (thr) at (1.5, -3) {};
       \node[dot] (cou3) at (2.5, -3.5) {};
       \node[dot] (cou2) at (1.5, -4.5) {};
       \draw (sta1) to[out=270, in=90] (zet.north west);
       \draw (sta2) to[out=270, in=90] (one.north west);
       \draw (sta3) to[out=270, in=90] (one.north east);
       \draw (sta4) to[out=270, in=0] (cou3.center);
       \draw (sta5) to[out=270, in=0] (cou2.center);
       \draw (one.south west) to[out=270, in=90] (zet.north east);
       \draw (one.south east) to[out=315, in=45] node[right]{$\nabla^*$} (thr.north east);
       \draw (zet.south west) to[out=270, in=90] (end1);
       \draw (zet.south east) to[out=270, in=90] (thr.north west);
       \draw (thr.south west) to[out=270, in=180] (cou2.center);
       \draw (thr.south east) to[out=270, in=180] (cou3.center);
     \end{tikzpicture}
   \end{align*}
   Thus \cref{EqNuAssoc2} holds, and so does \cref{EqNuAssoc}.
\end{proof}

By the universal property of $\ffunctor{M}(\group{S})^{\monad{T}}$
(\cref{propEMModule}) applied to the left $\monad{T}$-module $(L, \cone{L}_S, \nu)$, we have:
\begin{proposition}\label{propDefA}
 There is a unique functor ${A \colon L \to
 \ffunctor{M}(\group{S})^{\monad{T}}}$ fitting in the diagram
 \[
   \begin{tikzcd}
     L \arrow[rr, "A", dotted] \arrow[rd, "\cone{L}_{\group{S}}", swap] & &
     \ffunctor{M}(S)^{\monad{T}} \arrow[ld, "U"] \\
     & \ffunctor{M}(S) &
   \end{tikzcd}
 \]
 and satisfying $\bar{\mu} A = \nu$.
\end{proposition}

\subsection{The comparison functor $\ffunctor{M}(\group{S})^{\monad{T}} \to L$}

In the reverse direction, giving a functor $B \colon \ffunctor{M}(\group{S})^{\monad{T}}
\to L$ is equivalent (\cref{defBilimit}) to giving a pseudonatural transformation
$\phi \colon \Delta\ffunctor{M}(\group{S})^{\monad{T}} \Rightarrow \ffunctor{M}
\circ \ffunctor{U}$ between 2-functors $\cattransportext{\group{G}}{\group{S}} \to \catcat$.
Define the components of $\phi$ as follows:
\begin{itemize}
\item for any object $(\groupoid{P}, j_{\groupoid{P}})$ of
  $\cattransportext{\group{G}}{\group{S}}$,
  \[
    \begin{tikzcd}
      \phi_{(\groupoid{P}, j_{\groupoid{P}})} \colon
      \ffunctor{M}(\group{S})^{\monad{T}} \arrow[r, "U"] &
      \ffunctor{M}(\group{S}) \arrow[r, "j_{\groupoid{P}}^*"] &
      \ffunctor{M}(\groupoid{P})
    \end{tikzcd}
  \]
\item for any $1$-morphism $(a, \alpha) \colon (\groupoid{P}, j_{\groupoid{P}}) \to
  (\groupoid{Q}, j_{\groupoid{Q}})$,
  that is a $1$-morphism $a \colon \groupoid{P} \to \groupoid{Q}$ and a
  $2$-morphism $\alpha \colon ij_{\groupoid{Q}}a \Rightarrow ij_{\groupoid{P}}$
  in $\catgpdf$,
  \begin{equation}\label{defPhiMorph}
    \phi_{(a, \alpha)} \spacedequal
    \begin{tikzpicture}[baseline=(current bounding box.center)]
      \node[tip] (sta1) at (0, 0) {$U$}; \node[tip] (sta2) at (1.5, 0)
      {$j_{\groupoid{Q}}^*$}; \node[tip] (sta3) at (2, 0) {$a^*$};

      \node[dot] (uni) at (0.75, -0.5) {}; \node[naturaltr] (alp) at (1.5, -1)
      {$\alpha^*$}; \node[naturaltr] (mu) at (0.5, -2) {$\bar{\mu}$};

      \node[tip] (end1) at (0.5, -3) {$U$}; \node[tip] (end2) at (1.5, -3)
      {$j_{\groupoid{P}}^*$};

      \draw (sta1) to[out=270, in=90] (mu.north west); \draw (sta2) to[out=270,
      in=90] (alp.north); \draw (sta3) to[out=270, in=90] (alp.north east);
      \draw (uni.center) to[out=180, in=90] (mu.north); \draw (uni.center)
      to[out=0, in=90] (alp.north west); \draw (alp.south west) to[out=270,
      in=90] (mu.north east); \draw (alp.south) to[out=270, in=90] (end2); \draw
      (mu.south) to[out=270, in=90] (end1);
    \end{tikzpicture}
  \end{equation}
\end{itemize}
\begin{proposition}
 The natural transformation $\phi_{(a, \alpha)}$ is invertible.
\end{proposition}
\begin{proof}
 Consider
 \[
   \psi_{(a, \alpha)} \spacedequal
   \begin{tikzpicture}[baseline=(current bounding box.center)]
    \node[tip] (sta1) at (0, 0) {$U$};
    \node[tip] (sta2) at (1.5, 0) {$j_{\groupoid{P}}^*$};

    \node[dot] (uni) at (0.75, -0.5) {};
    \node[naturaltr] (alp) at (1.5, -1) {$(\alpha^*)^{-1}$};
    \node[naturaltr] (mu) at (0.5, -2) {$\bar{\mu}$};

    \node[tip] (end1) at (0.5, -3) {$U$};
    \node[tip] (end2) at (1.5, -3) {$j_{\groupoid{Q}}^*$};
    \node[tip] (end3) at (2, -3) {$a^*$};

    \draw (sta1) to[out=270, in=90] (mu.north west);
    \draw (sta2) to[out=270, in=90] (alp.north);
    \draw (uni.center) to[out=180, in=90] (mu.north);
    \draw (uni.center) to[out=0, in=90] (alp.north west);
    \draw (alp.south west) to[out=270, in=90] (mu.north east);
    \draw (alp.south) to[out=270, in=90] (end2);
    \draw (alp.south east) to[out=270, in=90] (end3);
    \draw (mu.south) to[out=270, in=90] (end1);
  \end{tikzpicture}
\]
We can check, using the associativity and the unitality of the action
$\bar{\mu}$, and the definition of the multiplication $\mu$ of $\monad{T}$, that $\psi_{(a, \alpha)} \circ \phi_{(a, \alpha)} = \id$ and
$\phi_{(a, \alpha)} \circ \psi_{(a, \alpha)} = \id$.
\end{proof}

\begin{proposition}
 The family of functors and natural transformations $\phi$ is a
 pseudonatural transformation ${\phi \colon \Delta\ffunctor{M}(\group{S})^{\monad{T}} \Rightarrow \ffunctor{M}
\circ \ffunctor{U}}$.
\end{proposition}
\begin{proof}
  This is done by the following straightforward computations.
  \[
    \phi_{(\id,\id)}
    \spacedequal
    \begin{tikzpicture}[baseline=(current bounding box.center)]
      \node[tip] (sta1) at (0, 0) {$U$};
      \node[tip] (sta2) at (1.5, 0) {$j_{\groupoid{P}}^*$};
      \node[dot] (uni) at (0.75, -0.5) {};
      \node[naturaltr] (alp) at (1.5, -1) {$\id^*$};
      \node[naturaltr] (mu) at (0.5, -2) {$\bar{\mu}$};
      \node[tip] (end1) at (0.5, -3) {U};
      \node[tip] (end2) at (1.5, -3) {$j_{\groupoid{P}}^*$};
      \draw (sta1) to[out=270, in=90] (mu.north west);
      \draw (sta2) to[out=270, in=90] (alp.north);
      \draw (uni.center) to[out=180, in=90] (mu.north);
      \draw (uni.center) to[out=0, in=90] (alp.north west);
      \draw (alp.south west) to[out=270, in=90] (mu.north east);
      \draw (alp.south) to[out=270, in=90] (end2);
      \draw (mu.south) to[out=270, in=90] (end1);
    \end{tikzpicture}
    \spacedequal
    \begin{tikzpicture}[baseline=(current bounding box.center)]
      \node[tip] (sta1) at (0, 0) {$U$};
      \node[tip] (sta2) at (1.5, 0) {$j_{\groupoid{P}}^*$};
      \node[dot] (uni) at (0.75, -0.5) {};
      \node[naturaltr] (mu) at (0.5, -2) {$\bar{\mu}$};
      \node[tip] (end1) at (0.5, -3) {$U$};
      \node[tip] (end2) at (1.5, -3) {$j_{\groupoid{P}}^*$};
      \draw (sta1) to[out=270, in=90] (mu.north west);
      \draw (sta2) to[out=270, in=90] (end2);
      \draw (uni.center) to[out=180, in=90] (mu.north);
      \draw (uni.center) to[out=0, in=90] (mu.north east);
      \draw (mu.south) to[out=270, in=90] (end1);
    \end{tikzpicture}
    \spacedequal
    \begin{tikzpicture}[baseline=(current bounding box.center)]
      \node[tip] (sta1) at (0, 0) {$U$};
      \node[tip] (sta2) at (1, 0) {$j_{\groupoid{P}}^*$};
      \node[tip] (end1) at (0, -3) {$U$};
      \node[tip] (end2) at (1, -3) {$j_{\groupoid{P}}^*$};
      \draw (sta1) to[out=270, in=90] (end1);
      \draw (sta2) to[out=270, in=90] (end2);
    \end{tikzpicture}
    \spacedequal \id_{j_{\groupoid{P}^*U}}
  \]
  For any composable morphisms ${(a, \alpha) \colon (\groupoid{P}, j_{\groupoid{P}})
  \to (\groupoid{Q}, j_{\groupoid{Q}})}$ and ${(b, \beta) \colon (\groupoid{Q},
  j_{\groupoid{Q}}) \to (\groupoid{R}, j_{\groupoid{R}})}$, using the
associativity of the action $\bar{\mu}$ and the definition of the multiplication
$\mu$ of the monad $\monad{T}$:
  \begin{align*}
    \phi_{(a,\alpha)} \circ a^*\phi_{(b, \beta)}
    &\spacedequal
     \begin{tikzpicture}[baseline=(current bounding box.center)]
      \node[tip] (sta1) at (0, 0) {$U$};
      \node[tip] (sta2) at (1.5, 0) {$j_{\groupoid{R}}^*$};
      \node[tip] (sta3) at (2, 0) {$b^*$};
      \node[tip] (sta4) at (2.5, 0) {$a^*$};
      \node[dot] (uni1) at (0.75, -0.5) {};
      \node[naturaltr] (bet) at (1.5, -1) {$\beta^*$};
      \node[naturaltr] (mu1) at (0.5, -2) {$\bar{\mu}$};
      \node[dot] (uni2) at (1.25, -2.5) {};
      \node[naturaltr] (alp) at (2, -3) {$\alpha^*$};
      \node[naturaltr] (mu2) at (1, -4) {$\bar{\mu}$};
      \node[tip] (end1) at (1, -5) {$U$};
      \node[tip] (end2) at (2, -5) {$j_{\groupoid{P}}^*$};
      \draw (sta1) to[out=270, in=90] (mu1.north west);
      \draw (sta2) to[out=270, in=90] (bet.north);
      \draw (sta3) to[out=270, in=90] (bet.north east);
      \draw (sta4) to[out=270, in=90] (alp.north east);
      \draw (uni1.center) to[out=180, in=90] (mu1.north);
      \draw (uni1.center) to[out=0, in=90] (bet.north west);
      \draw (bet.south west) to[out=270, in=90] (mu1.north east);
      \draw (bet.south) to[out=270, in=90] (alp.north);
      \draw (mu1.south) to[out=270, in=90] (mu2.north west);
      \draw (uni2.center) to[out=180, in=90] (mu2.north);
      \draw (uni2.center) to[out=0, in=90] (alp.north west);
      \draw (alp.south west) to[out=270, in=90] (mu2.north east);
      \draw (alp.south) to[out=270, in=90] (end2);
      \draw (mu2.south) to[out=270, in=90] (end1);
    \end{tikzpicture}
    \spacedequal
     \begin{tikzpicture}[baseline=(current bounding box.center)]
      \node[tip] (sta1) at (0, 0) {$U$};
      \node[tip] (sta2) at (1.5, 0) {$j_{\groupoid{R}}^*$};
      \node[tip] (sta3) at (2, 0) {$b^*$};
      \node[tip] (sta4) at (2.5, 0) {$a^*$};
      \node[dot] (uni1) at (0.75, -0.5) {};
      \node[naturaltr] (bet) at (1.5, -1) {$\beta^*$};
      \node[naturaltr] (mu1) at (1, -3.25) {$\mu$};
      \node[dot] (uni2) at (1.25, -1.75) {};
      \node[naturaltr] (alp) at (2, -2.25) {$\alpha^*$};
      \node[naturaltr] (mu2) at (0.5, -4) {$\bar{\mu}$};
      \node[tip] (end1) at (0.5, -5) {$U$};
      \node[tip] (end2) at (2, -5) {$j_{\groupoid{P}}^*$};
      \draw (sta1) to[out=270, in=90] (mu2.north west);
      \draw (sta2) to[out=270, in=90] (bet.north);
      \draw (sta3) to[out=270, in=90] (bet.north east);
      \draw (sta4) to[out=270, in=90] (alp.north east);
      \draw (uni1.center) to[out=180, in=90] (mu1.north west);
      \draw (uni1.center) to[out=0, in=90] (bet.north west);
      \draw (bet.south west) to[out=270, in=90] (mu1.135);
      \draw (bet.south) to[out=270, in=90] (alp.north);
      \draw (uni2.center) to[out=180, in=90] (mu1.north);
      \draw (uni2.center) to[out=0, in=90] (alp.north west);
      \draw (alp.south west) to[out=270, in=90] (mu1.north east);
      \draw (alp.south) to[out=270, in=90] (end2);
      \draw (mu1.south west) to[out=270, in=90] (mu2.north);
      \draw (mu1.south) to[out=270, in=90] (mu2.north east);
      \draw (mu2.south) to[out=270, in=90] (end1);
    \end{tikzpicture} \\
  &\spacedequal \begin{tikzpicture}[baseline=(current bounding box.center)]
      \node[tip] (sta1) at (0, 0) {$U$};
      \node[tip] (sta2) at (1.5, 0) {$j_{\groupoid{R}}^*$};
      \node[tip] (sta3) at (2, 0) {$b^*$};
      \node[tip] (sta4) at (2.5, 0) {$a^*$};
      \node[dot] (uni1) at (0.75, -0.5) {};
      \node[naturaltr] (bet) at (1.5, -1) {$\beta^*$};
      \node[dot] (cou1) at (1, -2.5) {};
      \node[dot] (uni2) at (1.25, -1.75) {};
      \node[naturaltr] (alp) at (2, -2.25) {$\alpha^*$};
      \node[naturaltr] (mu2) at (0.5, -3.5) {$\bar{\mu}$};
      \node[tip] (end1) at (0.5, -4.5) {$U$};
      \node[tip] (end2) at (2, -4.5) {$j_{\groupoid{P}}^*$};
      \draw (sta1) to[out=270, in=90] (mu2.north west);
      \draw (sta2) to[out=270, in=90] (bet.north);
      \draw (sta3) to[out=270, in=90] (bet.north east);
      \draw (sta4) to[out=270, in=90] (alp.north east);
      \draw (uni1.center) to[out=180, in=90] (mu2.north);
      \draw (uni1.center) to[out=0, in=90] (bet.north west);
      \draw (bet.south west) to[out=270, in=180] (cou1.center);
      \draw (bet.south) to[out=270, in=90] (alp.north);
      \draw (uni2.center) to[out=180, in=0] (cou1.center);
      \draw (uni2.center) to[out=0, in=90] (alp.north west);
      \draw (alp.south west) to[out=270, in=90] (mu2.north east);
      \draw (alp.south) to[out=270, in=90] (end2);
      \draw (mu2.south) to[out=270, in=90] (end1);
    \end{tikzpicture}
    \spacedequal \begin{tikzpicture}[baseline=(current bounding box.center)]
      \node[tip] (sta1) at (0, 0) {$U$};
      \node[tip] (sta2) at (1.25, 0) {$j_{\groupoid{R}}^*$};
      \node[tip] (sta3) at (1.75, 0) {$b^*$};
      \node[tip] (sta4) at (2.25, 0) {$a^*$};
      \node[dot] (uni1) at (0.75, -0.5) {};
      \node[naturaltr] (alp) at (1.5, -2) {$\alpha^* \circ a^*\beta^*$};
      \node[naturaltr] (mu2) at (0.5, -3.5) {$\bar{\mu}$};
      \node[tip] (end1) at (0.5, -4.5) {$U$};
      \node[tip] (end2) at (1.5, -4.5) {$j_{\groupoid{P}}^*$};
      \draw (sta1) to[out=270, in=90] (mu2.north west);
      \draw (sta2) to[out=270, in=90] (alp.110);
      \draw (sta3) to[out=270, in=90] (alp.60);
      \draw (sta4) to[out=270, in=90] (alp.25);
      \draw (uni1.center) to[out=180, in=90] (mu2.north);
      \draw (uni1.center) to[out=0, in=90] (alp.140);
      \draw (alp.south west) to[out=250, in=90] (mu2.north east);
      \draw (alp.south) to[out=270, in=90] (end2);
      \draw (mu2.south) to[out=270, in=90] (end1);
    \end{tikzpicture} \\
    &\spacedequal \phi_{(b, \beta) \circ (a, \alpha)}
  \end{align*}
  For any $2$-morphism $\zeta \colon (a, \alpha) \Rightarrow (b, \beta)$ between
  parallel morphisms $(\groupoid{P}, j_{\groupoid{P}}) \to (\groupoid{Q},
  j_{\groupoid{Q}})$ in $\cattransportext{G}{S}$, we compute:
  \begin{align*}
    \phi_{(b, \beta)} \circ \phi_{(\groupoid{Q},j_{\groupoid{Q}})}(\zeta^*) =
    \begin{tikzpicture}[baseline=(current bounding box.center)]
      \node[tip] (sta1) at (0, 0) {$U$};
      \node[tip] (sta2) at (1.5, 0) {$j_{\groupoid{Q}}^*$};
      \node[tip] (sta3) at (2.5, 0) {$a^*$};
      \node[naturaltr] (zet) at (2.5,-1) {$\zeta^*$};
      \node[dot] (uni) at (0.75, -1.5) {};
      \node[naturaltr] (alp) at (1.5, -2) {$\beta^*$};
      \node[naturaltr] (mu) at (0.5, -3) {$\bar{\mu}$};
      \node[tip] (end1) at (0.5, -4) {$U$};
      \node[tip] (end2) at (1.5, -4) {$j_{\groupoid{P}}^*$};
      \draw (sta1) to[out=270, in=90] (mu.north west);
      \draw (sta2) to[out=270, in=90] (alp.north);
      \draw (sta3) to[out=270, in=90] (zet.north);
      \draw (zet.south west) to[out=270, in=90] (alp.north east);
      \draw (uni.center) to[out=180, in=90] (mu.north);
      \draw (uni.center) to[out=0, in=90] (alp.north west);
      \draw (alp.south west) to[out=270, in=90] (mu.north east);
      \draw (alp.south) to[out=270, in=90] (end2);
      \draw (mu.south) to[out=270, in=90] (end1);
    \end{tikzpicture} =
    \begin{tikzpicture}[baseline=(current bounding box.center)]
      \node[tip] (sta1) at (0, 0) {$U$};
      \node[tip] (sta2) at (1.5, 0) {$j_{\groupoid{Q}}^*$};
      \node[tip] (sta3) at (2, 0) {$a^*$};
      \node[dot] (uni) at (0.75, -0.5) {};
      \node[naturaltr] (alp) at (1.5, -1) {$\alpha^*$};
      \node[naturaltr] (mu) at (0.5, -2) {$\bar{\mu}$};
      \node[tip] (end1) at (0.5, -3) {$U$};
      \node[tip] (end2) at (1.5, -3) {$j_{\groupoid{P}}^*$};
      \draw (sta1) to[out=270, in=90] (mu.north west);
      \draw (sta2) to[out=270, in=90] (alp.north);
      \draw (sta3) to[out=270, in=90] (alp.north east);
      \draw (uni.center) to[out=180, in=90] (mu.north);
      \draw (uni.center) to[out=0, in=90] (alp.north west);
      \draw (alp.south west) to[out=270, in=90] (mu.north east);
      \draw (alp.south) to[out=270, in=90] (end2);
      \draw (mu.south) to[out=270, in=90] (end1);
    \end{tikzpicture} \quad = \phi_{(a,\alpha)}
  \end{align*}
  These equalities show the pseudonaturality of $\phi$.
\end{proof}

Hence, by the universal property of the bilimit:
\begin{proposition}\label{propDefB}
  There exists a unique (up-to 2-isomorphism) functor
  \[
    B \colon \ffunctor{M}(\group{S})^{\monad{T}} \to L
  \] and a modification $m \colon \phi \cong \cone{L} \cdot B$.
\end{proposition}

\subsection{The functors $A$ and $B$ are mutually pseudoinverse}

\begin{proposition}\label{propABId}
  The composite $AB \colon \ffunctor{M}(\group{S})^{\monad{T}} \to
 \ffunctor{M}(\group{S})^{\monad{T}}$ is isomorphic to the identity
 functor.
\end{proposition}
\begin{proof}
 The universal property (\cref{remarkEMBiterm}) of
 $\ffunctor{M}(\group{S})^{\monad{T}}$ guarantees that
 $(\id_{\ffunctor{M}(S)^{\monad{T}}}, \id_U)$ is the unique $1$-endomorphism of the
 left $\monad{T}$-module $(\ffunctor{M}(S)^{\monad{T}}, U, \bar{\mu})$, up to a
 $2$-isomorphism. Hence, the existence of a $2$-isomorphism $\xi \colon U \xRightarrow{\sim} UAB$ such that $(AB, \xi) :
 (\ffunctor{M}(S)^{\monad{T}},U,\bar{\mu}) \to
 (\ffunctor{M}(S)^{\monad{T}},U,\bar{\mu})$ is a $1$-morphism of left
 $\monad{T}$-modules (\cref{defModule2Cat}) implies that
 $\id_{\ffunctor{M}(S)^{\monad{T}}}$ and $AB$ are isomorphic. Let us find such a $\xi$.

 First, as
 \[
   \phi_{\group{S}} = j_{\group{S}}^*U = U
 \]
 we can take the $2$-isomorphism
 \[
   m_{\group{S}} \colon UAB = \cone{L}_{\group{S}}B \cong \phi_S = U
 \]
 given by \cref{propDefB} as our $2$-isomorphism $\xi$.
 Moreover, since $m \colon \cone{L}B \to \phi$ is a modification,
 \begin{align*}
   &m_{\group{S}} \circ \bar{\mu}AB \overset{\eqref{propDefA}}{=} m_{\group{S}} \circ \nu B
   \xspacedequal{\eqref{defNu}}
   \begin{tikzpicture}[baseline=(current bounding box.center)]
     \node[tip] (sta4) at (-1, 0) {$B$};
     \node[tip] (sta1) at (0,0) {$\cone{L}_{\group{S}}$};
    \node[tip] (sta2) at (1,0) {$i_!$};
    \node[tip] (sta3) at (2,0) {$i^*$};
    \node[naturaltr] (lam1) at (1.5, -1) {$\lambda_!^{-1}$};
    \node[naturaltr] (lam2) at (0.5,-2) {$\cone{L}_\lambda$};
    \node[dot] (cou2) at (1.5,-2.5) {};
    \node[naturaltr] (m) at (-0.5, -3) {$m_{\group{S}}$} ;
    \node[tip] (end) at (-0.5, -4) {$U$};
    \draw (sta4)            to[out=270, in=90]  (m.north west);
    \draw (sta1)            to[out=270, in=90]  (lam2.north west);
    \draw (sta2)            to[out=270, in=90]  (lam1.north west);
    \draw (sta3)            to[out=270, in=90]  (lam1.north east);
    \draw (lam1.south west) to[out=270, in=90]  (lam2.north east);
    \draw (lam1.south east) to[out=270, in=0]   (cou2.center);
    \draw (lam2.south west) to[out=270, in=90]  (m.north east);
    \draw (lam2.south east) to[out=270, in=180] (cou2.center);
    \draw (m.south)         to[out=270, in=90]  (end);
  \end{tikzpicture} \\
  & \xspacedequal{\eqref{propDefB}}
  \begin{tikzpicture}[baseline=(current bounding box.center)]
    \node[tip] (sta4) at (-1, 0) {$B$};
    \node[tip] (sta1) at (0,0) {$\cone{L}_{\group{S}}$};
    \node[tip] (sta2) at (1,0) {$i_!$};
    \node[tip] (sta3) at (2,0) {$i^*$};
    \node[naturaltr] (lam1) at (1.5, -1) {$\lambda_!^{-1}$};
    \node[naturaltr] (m) at (-0.5, -1) {$m_{\group{S}}$};
    \node[naturaltr] (lam2) at (0.5,-2) {$\phi_{(\id,\lambda)}$};
    \node[dot] (cou2) at (1.5,-2.5) {};
    \node[tip] (end1) at (0,-3) {$U$};
    \draw (sta4)            to[out=270, in=90]  (m.north west);
    \draw (sta1)            to[out=270, in=90]  (m.north east);
    \draw (sta2)            to[out=270, in=90]  (lam1.north west);
    \draw (sta3)            to[out=270, in=90]  (lam1.north east);
    \draw (lam1.south west) to[out=270, in=90]  (lam2.north east);
    \draw (lam1.south east) to[out=270, in=0]   (cou2.center);
    \draw (m.south)         to[out=270, in=90]  (lam2.north west);
    \draw (lam2.south west) to[out=270, in=90]  (end1);
    \draw (lam2.south east) to[out=270, in=180] (cou2.center);
  \end{tikzpicture}
  \xspacedequal{\eqref{defPhiMorph}}
  \begin{tikzpicture}[baseline=(current bounding box.center)]
    \node[tip] (sta4) at (-1, 0) {$B$};
    \node[tip] (sta1) at (0,0) {$\cone{L}_{\group{S}}$};
    \node[tip] (sta2) at (2,0) {$i_!$};
    \node[tip] (sta3) at (3,0) {$i^*$};
    \node[naturaltr] (lam1) at (2.5, -1) {$\lambda_!^{-1}$};
    \node[naturaltr] (m) at (-0.5, -1) {$m_{\group{S}}$};
    \node[dot] (uni2) at (0.75,-2) {};
    \node[naturaltr] (lam2) at (1.5, -2.5) {$\lambda^*$};
    \node[dot] (cou2) at (2.5,-3) {};
    \node[naturaltr] (mu) at (0,-3.5) {$\bar{\mu}$};
    \node[tip] (end1) at (0,-4.5) {$U$};
    \draw (sta4)            to[out=270, in=90]  (m.north west);
    \draw (sta1)            to[out=270, in=90]  (m.north east);
    \draw (sta2)            to[out=270, in=90]  (lam1.north west);
    \draw (sta3)            to[out=270, in=90]  (lam1.north east);
    \draw (lam1.south west) to[out=270, in=90]  (lam2.north east);
    \draw (lam1.south east) to[out=270, in=0]   (cou2.center);
    \draw (m.south)         to[out=270, in=90]  (mu.north west);
    \draw (uni2.center)     to[out=180, in=90]  (mu.55);
    \draw (uni2.center)     to[out=0,   in=90]  (lam2.north west);
    \draw (lam2.south west) to[out=270, in=90]  (mu.35);
    \draw (lam2.south east) to[out=270, in=180] (cou2.center);
    \draw (mu.270)          to[out=270, in=90]  (end1);
  \end{tikzpicture} \\
  &\xspacedequal{\eqref{eqDefMate}}
  \begin{tikzpicture}[baseline=(current bounding box.center)]
    \node[tip] (sta4) at (-1, 0) {$B$};
    \node[tip] (sta1) at (0,0) {$\cone{L}_{\group{S}}$};
    \node[tip] (sta2) at (1,0) {$i_!$};
    \node[tip] (sta3) at (2,0) {$i^*$};
    \node[naturaltr] (m) at (-0.5, -1) {$m_{\group{S}}$};
    \node[naturaltr] (mu) at (1,-2) {$\bar{\mu}$};
    \node[tip] (end1) at (1,-3) {$U$};
    \draw (sta4)            to[out=270, in=90] (m.north west);
    \draw (sta1)            to[out=270, in=90] (m.north east);
    \draw (m.south)         to[out=270, in=90] (mu.north west);
    \draw (sta2)            to[out=270, in=90] (mu.north);
    \draw (sta3)            to[out=270, in=90]  (mu.north east);
    \draw (mu.270)          to[out=270, in=90] (end1);
  \end{tikzpicture} \spacedequal \bar{\mu} \circ \monad{T}m_{\group{S}}
 \end{align*}
 This equality shows that $(AB, m_S) \colon (\ffunctor{M}(S)^{\monad{T}}, U, \bar{\mu}) \to
 (\ffunctor{M}(S)^{\monad{T}}, U, \bar{\mu})$ is a $1$-morphism of left
 $\monad{T}$-modules. Since such a $1$-morphism is unique up to $2$-isomorphism, $AB$ is isomorphic to $\id_{\ffunctor{M}(\group{S})^{\monad{T}}}$.
\end{proof}

\begin{proposition}\label{propBAId}
  The composite $BA : L \to L$ is isomorphic to the identity functor.
\end{proposition}
\begin{proof}
 Similarly we will show that the induced cone $\tilde{\cone{L}} := \cone{L} \cdot
 (BA)$ is isomorphic to the standard cone $\cone{L}$, that is, that there
 is an invertible modification $n$ between them.
 On any object $(\groupoid{P}, j_{\groupoid{P}})$ define $n_{\groupoid{P}}$ as
 \[
   n_{\groupoid{P}} \colon \tilde{\cone{L}}_{\groupoid{P}} =
   \cone{L}_{\groupoid{P}}BA \xrightarrow{m_{\groupoid{P}}A} \phi_{\groupoid{P}} A = j_{\groupoid{P}}^*UA =
   j_{\groupoid{P}}^*\cone{L}_S = \cone{L}_{\groupoid{P}}
 \]
 The family of morphisms $(n_{\groupoid{P}})$ is an invertible modification
 between the cones $\tilde{\cone{L}}$ and $\phi A$. Thus, it suffices to check
 that $\phi A$ and $\cone{L}$ are indeed the same cones. Their components
 at any object $(\groupoid{P}, j_{\groupoid{P}})$ coincide; we still have to
 check that the components at any morphism $(a, \alpha)$ are equal.

 At the morphism $(\id, \lambda) : (\ppback{i}{i}, \ppbackr{i}{i}) \rightarrow
 (\ppback{i}{i}, \ppbackl{i}{i})$ of $\cattransportext{G}{S}$, we have:
 \begin{align*}
  &\phi_{(\id, \lambda)}A \spacedequal
  \begin{tikzpicture}[baseline=(current bounding box.center)]
    \node[tip] (sta1) at (0, 0) {$A$};
    \node[tip] (sta2) at (1, 0) {$U$};
    \node[tip] (sta3) at (2.5, 0) {$\ppbackl{i}{i}^*$};
    \node[dot] (uni) at (1.75, -0.5) {};
    \node[naturaltr] (alp) at (2.5, -1) {$\lambda^*$};
    \node[naturaltr] (mu) at (1.5, -2) {$\bar{\mu}$};
    \node[tip] (end1) at (0, -3) {$A$};
    \node[tip] (end2) at (1.5, -3) {$U$};
    \node[tip] (end3) at (2.5, -3) {$\ppbackr{i}{i}^*$};
    \draw (sta1) to[out=270, in=90]  (end1);
    \draw (sta2) to[out=270, in=90] (mu.north west);
    \draw (sta3) to[out=270,in=90] (alp.north);
    \draw (uni.center) to[out=180, in=90] (mu.north);
    \draw (uni.center) to[out=0, in=90] (alp.north west);
    \draw (alp.south west) to[out=270, in=90] (mu.north east);
    \draw (alp.south) to[out=270, in=90] (end3);
    \draw (mu.south) to[out=270, in=90] (end2);
  \end{tikzpicture}
  \xspacedequal{\eqref{propDefA}}
  \begin{tikzpicture}[baseline=(current bounding box.center)]
    \node[tip] (sta1) at (0.5, 0) {$\cone{L}_{\group{S}}$};
    \node[tip] (sta2) at (2.5, 0) {$\ppbackl{i}{i}^*$};
    \node[dot] (uni1) at (1.25, -0.5) {};
    \node[naturaltr] (alp) at (2, -1) {$\lambda^*$};
    \node[naturaltr] (nu) at (1, -2) {$\nu$};
    \node[tip] (end1) at (1, -3) {$\cone{L}_{\group{S}}$};
    \node[tip] (end2) at (2.5, -3) {$\ppbackr{i}{i}^*$};
    \draw (sta1) to[out=270, in=90]  (nu.north west);
    \draw (sta2) to[out=270, in=90] (alp.north east);
    \draw (uni1.center) to[out=180, in=90] (nu.north);
    \draw (uni1.center) to[out=0, in=90] (alp.north west);
    \draw (alp.south west) to[out=270, in=90] (nu.north east);
    \draw (alp.south east) to[out=270, in=90] (end2);
    \draw (nu.south) to[out=270, in=90] (end1);
  \end{tikzpicture} \\
  & \xspacedequal{\eqref{defNu}}
  \begin{tikzpicture}[baseline=(current bounding box.center)]
    \node[tip] (sta1) at (-0.5, 0) {$\cone{L}_{\group{S}}$};
    \node[tip] (sta2) at (2.5, 0) {$\ppbackl{i}{i}^*$};
    \node[dot] (uni1) at (1, -0.5) {};
    \node[naturaltr] (alp) at (2, -1) {$\lambda^*$};
    \node[naturaltr] (lam1) at (1, -2) {$\lambda_!^{-1}$};
    \node[naturaltr] (lam2) at (0, -3) {$\cone{L}_{\lambda}$};
    \node[dot] (cou2) at (1, -3.5) {};
    \node[tip] (end1) at (0, -4) {$\cone{L}_{\group{S}}$};
    \node[tip] (end2) at (2.5, -4) {$\ppbackr{i}{i}^*$};
    \draw (sta1) to[out=270, in=90]  (lam2.north west);
    \draw (sta2) to[out=270, in=90] (alp.north east);
    \draw (uni1.center) to[out=180, in=90] (lam1.north west);
    \draw (uni1.center) to[out=0, in=90] (alp.north west);
    \draw (alp.south west) to[out=270, in=90] (lam1.north east);
    \draw (alp.south east) to[out=270, in=90] (end2);
    \draw (lam1.south west) to[out=270, in=90] (lam2.north east);
    \draw (lam1.south east) to[out=270, in=0] (cou2.center);
    \draw (lam2.south) to[out=270, in=90] (end1);
    \draw (lam2.south east) to[out=270, in=180] (cou2.center);
  \end{tikzpicture}
  \xspacedequal{\eqref{eqDefMate}}
  \begin{tikzpicture}[baseline=(current bounding box.center)]
    \node[tip] (sta1) at (-0.5, 0) {$\cone{L}_{\group{S}}$};
    \node[tip] (sta2) at (2.5, 0) {$\ppbackl{i}{i}^*$};
    \node[dot] (uni1) at (1.25, -0.5) {};
    \node[dot] (uni3) at (3.5, -0.5) {};
    \node[naturaltr] (alp) at (2, -1) {$\lambda^*$};
    \node[dot] (cou3) at (2.75, -1.5) {};
    \node[naturaltr] (lam1) at (1, -2) {$\lambda_!^{-1}$};
    \node[naturaltr] (lam2) at (0, -3) {$\cone{L}_{\lambda}$};
    \node[dot] (cou2) at (1, -3.5) {};
    \node[tip] (end1) at (0, -4) {$\cone{L}_{\group{S}}$};
    \node[tip] (end2) at (4, -4) {$\ppbackr{i}{i}^*$};
    \draw (sta1) to[out=270, in=90]  (lam2.north west);
    \draw (sta2) to[out=270, in=90] (alp.north east);
    \draw (uni1.center) to[out=180, in=90] (lam1.north west);
    \draw (uni1.center) to[out=0, in=90] (alp.north west);
    \draw (uni3.center) to[out=0, in=90]  (end2);
    \draw (alp.south west) to[out=270, in=90] (lam1.north east);
    \draw (alp.south east) to[out=270, in=180] (cou3.center);
    \draw (cou3.center) to[out=0, in=180] (uni3.center);
    \draw (lam1.south west) to[out=270, in=90] (lam2.north east);
    \draw (lam1.south east) to[out=270, in=0] (cou2.center);
    \draw (lam2.south) to[out=270, in=90] (end1);
    \draw (lam2.south east) to[out=270, in=180] (cou2.center);
    %
  \end{tikzpicture} \\
  & \xspacedequal{\eqref{eqDefMate}}
  \begin{tikzpicture}[baseline=(current bounding box.center)]
    \node[tip] (sta1) at (0, 0) {$\cone{L}_{\group{S}}$};
    \node[tip] (sta2) at (1, 0) {$\ppbackl{i}{i}^*$};
    \node[dot] (uni2) at (1.75, -0.5) {};
    \node[naturaltr] (lam2) at (0.5, -1) {$\cone{L}_{\lambda}$};
    \node[dot] (cou2) at (1.25, -1.5) {};
    \node[tip] (end1) at (0.5, -2) {$\cone{L}_{\group{S}}$};
    \node[tip] (end2) at (2, -2) {$\ppbackr{i}{i}^*$};
    \draw (sta1) to[out=270, in=90] (lam2.north west);
    \draw (sta2) to[out=270, in=90] (lam2.north east);
    \draw (uni2.center) to[out=180, in=0] (cou2.center);
    \draw (uni2.center) to[out=0, in=90] (end2);
    \draw (lam2.south) to[out=270, in=90] (end1);
    \draw (lam2.south east) to[out=270, in=180] (cou2.center);
  \end{tikzpicture}
  \spacedequal \cone{L}_{\lambda}
 \end{align*}
For any morphism $(a, \alpha) : (\groupoid{P}, j_{\groupoid{P}}) \to
(\groupoid{Q}, j_{\groupoid{Q}})$ there is a (unique) morphism of groupoids
$\nabla_{\alpha} : \groupoid{P} \to \ppback{i}{i}$ such that
\begin{align*}
 \ppbackl{i}{i}\nabla_{\alpha} = j_{\groupoid{Q}}a \qquad \ppbackr{i}{i}\nabla_{\alpha} = j_{\groupoid{P}} \qquad \lambda\nabla_{\alpha} = \alpha
\end{align*}
Hence,
\begin{align*}
   \phi_{(a, \alpha)}A &\spacedequal
  \begin{tikzpicture}[baseline=(current bounding box.center)]
    \node[tip] (sta4) at (-1, 0) {$A$};
    \node[tip] (sta1) at (-0.25, 0) {$U$};
    \node[tip] (sta2) at (1.25, 0) {$j_{\groupoid{Q}}^*$};
    \node[tip] (sta3) at (1.75, 0) {$a^*$};
    \node[dot] (uni) at (0.5, -1.5) {};
    \node[naturaltr] (alp) at (1.25, -2) {$\alpha^*$};
    \node[naturaltr] (mu) at (0.25, -3) {$\bar{\mu}$};
    \node[tip] (end4) at (-1, -4) {$A$};
    \node[tip] (end1) at (0.25, -4) {$U$};
    \node[tip] (end2) at (1.25, -4) {$j_{\groupoid{P}}^*$};
    \draw (sta4) to[out=270, in=90]  (end4);
    \draw (sta1) to[out=270, in=90] (mu.north west);
    \draw (sta2) to[out=270, in=90] (alp.north);
    \draw (sta3) to[out=270, in=90] (alp.north east);
    \draw (uni.center) to[out=180, in=90] (mu.north);
    \draw (uni.center) to[out=0, in=90] (alp.north west);
    \draw (alp.south west) to[out=270, in=90] (mu.north east);
    \draw (alp.south) to[out=270, in=90] (end2);
    \draw (mu.south) to[out=270, in=90] (end1);
  \end{tikzpicture}
  \spacedequal
   \begin{tikzpicture}[baseline=(current bounding box.center)]
    \node[tip] (sta4) at (-1, 0) {$A$};
    \node[tip] (sta1) at (-0.25, 0) {$U$};
    \node[tip] (sta2) at (1.75, 0) {$j_{\groupoid{Q}}^*$};
    \node[tip] (sta3) at (2.25, 0) {$a^*$};
    \node[dot] (uni) at (0.5, -1.5) {};
    \node[iden] (id1) at (2, -1) {};
    \node[naturaltr] (alp) at (1.25, -2) {$\lambda^*$};
    \node[iden] (id2) at (2, -3) {};
    \node[naturaltr] (mu) at (0.25, -3) {$\bar{\mu}$};
    \node[tip] (end4) at (-1, -4) {$A$};
    \node[tip] (end1) at (0.25, -4) {$U$};
    \node[tip] (end2) at (2, -4) {$j_{\groupoid{P}}^*$};
    \draw (sta4) to[out=270, in=90]  (end4);
    \draw (sta1) to[out=270, in=90] (mu.north west);
    \draw (sta2) to[out=270, in=135] (id1.north west);
    \draw (sta3) to[out=270, in=45] (id1.north east);
    \draw (uni.center) to[out=180, in=90] (mu.north);
    \draw (uni.center) to[out=0, in=90] (alp.north west);
    \draw (id1.south west) to[out=225, in=90] (alp.north east);
    \draw (id1.south east) to[out=315, in=45] node[right]{$\nabla_{\alpha}^*$} (id2.north east);
    \draw (alp.south west) to[out=270, in=90] (mu.north east);
    \draw (alp.south east) to[out=270, in=135] (id2.north west);
    \draw (id2.270) to[out=270, in=90] (end2);
    \draw (mu.south) to[out=270, in=90] (end1);
  \end{tikzpicture} \\
  &\spacedequal
  \begin{tikzpicture}[baseline=(current bounding box.center)]
    \node[tip] (sta1) at (1.5,0)  {$a^*$};
    \node[tip] (sta2) at (1,0)  {$j_{\groupoid{Q}}^*$};
    \node[tip] (sta3) at (0,0)  {$\cone{L}_{\group{S}}$};
    \node[tip] (end2) at (1.25,-4) {$j_{\groupoid{P}}^*$};
     \node[tip] (end3) at (0,-4) {$\cone{L}_{\group{S}}$};
     \node[iden] (id1) at (1.25 ,-1) {};
     \node[naturaltr] (lam)  at (0.5  ,-2) {$\cone{L}_{\lambda}$};
     \node[iden] (id2) at (1.25 ,-3) {};
     \draw (sta1)    to[out=270, in=45 ] (id1.north east);
     \draw (sta2)    to[out=270, in=135] (id1.north west);
     \draw (sta3)    to[out=270, in=90] (lam.north west);
     \draw (id1.south east) to[out=315, in=45] node[right]{$\nabla_{\alpha}^*$} (id2.north east);
     \draw (id1.south west) to[out=225, in=90] (lam.north east);
     \draw (lam.south east) to[out=270, in=135] (id2.north west);
     \draw (id2.south) to[out=270, in=90 ] (end2);
     \draw (lam.south west) to[out=270, in=90 ] (end3);
   \end{tikzpicture}
  \spacedequal \cone{L}_{(a, \alpha)}
\end{align*}
Thus, $BA$ is isomorphic to the identity functor.
\end{proof}

Putting together \cref{propABId} and \cref{propBAId}, we get:
\begin{theorem}
 The categories $L$ and $\mackeyfunctor{M}(\group{S})^{\monad{T}}$ are equivalent.
\end{theorem}

\section{Applications} \label{secApp}

In this section, we describe two different ways of extracting the classical
Cartan-Eilenberg formula from \cref{theoremMackeyAs2Lim}.

\subsection{Extracting Hom-sets}

In this subsection, we will express the Hom-sets of a bilimit in $\catcat$ for a
\emph{strict} 2-functor as a certain limit, and show how it can be used to extract more concrete
information from the formula of \cref{theoremMackeyAs2Lim}.

We should first give an explicit description of pseudo bilimits in $\catcat$.
\begin{definition}
  Let $\bicatname{I}$ be a (small) $2$-category and $\ffunctor{D} \colon
  \bicatname{I} \to \catcat$ be a $2$-functor. Let
  $L_{\ffunctor{D}}$ be the category with:
  \begin{itemize}
  \item \emph{Objects:} The pairs of families $d := ((d_i)_i, (d_f)_f)$ such that:
    \begin{itemize}
    \item for each object $i \in \bicatname{I}$, $d_i$ is an object of
      $\ffunctor{D}i$
    \item for each $1$-morphism $f \colon i \to j$ in $\bicatname{I}$, $d_f$ is an isomorphism
      $\ffunctor{D}f(d_i) \xrightarrow{\sim} d_{j}$ in $\ffunctor{D}j$.
    \item for any $2$-morphism $\phi \colon f \Rightarrow g$ in $\bicatname{I}$,
      with $f,g \colon i \to j$,
      \[
        d_g \circ (\ffunctor{D}\phi)_{d_i} = d_f
      \]
    \item for any composable $1$-morphisms $f \colon i \to j$ and $g \colon j
      \to k$ in $\bicatname{I}$,
      \[
        d_{g \circ f} = d_g \circ (\ffunctor{D}g)(d_f)
      \]
    \end{itemize}
  \item \emph{Morphisms $((d_i), (d_f)) \to ((d_i'), (d_f'))$:} the families of
    morphisms
    \[
      (\delta_i \colon d_i \to d_i')_{i \in \Ob{\bicatname{I}}}
    \]
    such that, for any morphism $f \colon i \to j$ in $\bicatname{I}$, the
    following square commutes:
    \begin{equation}\label{eqMorphLD}
      \begin{tikzcd}
        \ffunctor{D}f(d_i) \arrow[r, "d_f"] \arrow[d, "\ffunctor{D}f(\delta_i)", swap]
        & d_j \arrow[d, "\delta_j"] \\
        \ffunctor{D}f(d_i') \arrow[r, "d_f'", swap] & d_j'
      \end{tikzcd}
    \end{equation}
  \item \emph{Composition} is induced componentwise by the compositions of the
    categories $(\ffunctor{D}i)_i$.
  \end{itemize}
  The category $L_{\ffunctor{D}}$ is endowed with a canonical cone $\phi$ over
  $\ffunctor{D}$ with components:
  \[
    \phi_i \colon \left\{
      \begin{array}{lll}
        L_{\ffunctor{D}} & \Rightarrow & \ffunctor{D}i \\
        ((d_i)_i, (d_f)_f) & \mapsto & d_i \\
        (f_i)_i & \mapsto & f_i
      \end{array}
    \right. \text{for any object } i \in \bicatname{I}
  \]
  \[
    (\phi_f)_{((d_i)_i, (d_f)_f)} = d_f \text{ for any morphism } f \colon i
    \to j \in \bicatname{I}
  \]
\end{definition}

\begin{proposition}
  The category $L_{\ffunctor{D}}$, with its canonical cone, is a bilimit of
  $\ffunctor{D}$:
  \[
    L_{\ffunctor{D}} \cong \llim_{\bicatname{I}}\ffunctor{D}
  \]
\end{proposition}
\begin{proof}
 The category $L_{\ffunctor{D}}$ is an explicit description of the category
 $[\Delta{}1, \ffunctor{D}]$ of the (pseudo)cones over $\ffunctor{D}$ with vertex $1$, the
 category with exactly one object and its identity morphism.

 In turn, $[\Delta{}1, \ffunctor{D}]$ is a model of the bilimit of $\ffunctor{D}$
 by the following classical equivalence, pseudonatural in $X \in \catcat$:
 \[
   \catcat(X, [\Delta{}1, \ffunctor{D}]) \cong \catcat(1, [\Delta{}X,
   \ffunctor{D}]) \cong [\Delta{}X, \ffunctor{D}]
 \]

\end{proof}

We now fix a $2$-diagram $\ffunctor{D} \colon \bicatname{I} \to \catcat$ and two
objects $d$, $d'$ of $L_{\ffunctor{D}}$. By definition:
\begin{equation}\label{eqExplicitHomLD}
  L_{\ffunctor{D}}(d, d') = \{ (\delta_i \colon d_i \to d_i')_i \: |
  \: \text{the squares \eqref{eqMorphLD} commute} \}
\end{equation}

We can define an associated $1$-diagram on the underlying $1$-category
$\bicatname{I}^{(1)}$ of $\bicatname{I}$, obtained by forgetting the $2$-morphisms,
as follows:
\begin{definition}
  The diagram $\ffunctor{D}_{d,d'} \colon \bicatname{I}^{(1)} \to \catset$ is:
  \[
    \ffunctor{D}_{d,d'} \colon \left\{
      \begin{array}{lll}
        i & \mapsto & \ffunctor{D}i(d_i, d_i') \\
        f : i \to j & \mapsto & \ffunctor{D}j(d_f^{-1}, d_f') \circ \ffunctor{D}f
      \end{array}
    \right.
  \]
\end{definition}

\begin{proposition}
  There is an isomorphism:
  \[
    L_{\ffunctor{D}}(d,d') \simeq \lim_{\bicatname{I}^{(1)}}\ffunctor{D}_{d,d'}
  \]
\end{proposition}
\begin{proof}
  The limit $\lim_{\bicatname{I}^{(1)}}\ffunctor{D}_{d,d'}$ has an explicit
  description, as a limit in $\catset$:
  \[
    \lim_{\bicatname{I}^{(1)}}\ffunctor{D}_{d,d'} = \{ (\delta_i \in
    \ffunctor{D}_{d,d'}(i))_i \: | \: \forall f : i \to j,
    \ffunctor{D}_{d,d'}(f)(\delta_i) = \delta_j \}
  \]
  Unfolding the definitions, we precisely get back the set of \eqref{eqExplicitHomLD}.
\end{proof}

\begin{remark}\label{remarkObjectTrivialAction}
  The expression of the diagram $\ffunctor{D}_{d,d'}$ simplifies when the
  objects d and d' are such that, for all $f \colon i \to j$, $d_f = d_f' = \id$:
  \[
    \ffunctor{D}_{d,d'} \colon \left\{
      \begin{array}{lll}
        i & \mapsto & \ffunctor{D}i(d_i, d_i') \\
        f : i \to j & \mapsto & \ffunctor{D}f
      \end{array}
    \right.
  \]
\end{remark}

\begin{example}\label{exampleCEHom}
  We can use these results to recover the Cartan-Eilenberg formula for the Tate
  cohomology $\cohomologyh{G}{\field}$ of a finite group $G$ over a field
  $\field$ of characteristic p. Indeed, for any finite group $G$, the Tate
  cohomology groups can be recovered as Hom-sets of $\catstmod{G}$:
  \[
    \catstmod{G}(\Omega^n\field,\field) = \ncohomologyh{n}{G}{\field}
  \]
  Fix a finite group $G$ and a $p$-Sylow $S$ of $G$. We consider the
  $2$-diagram:
  \[
    \ffunctor{D} = \catstmod{-} \circ \ffunctor{U} \colon \catop{\cattransportext{G}{S}}
    \to \catcat
  \]
  Note that, for any $n$, $\Omega^n\field \in \catstmod{G}$ can be seen as an object of
  $L_{\ffunctor{D}}$ satisfying the conditions of
  \cref{remarkObjectTrivialAction}. Hence we have:
  \begin{align*}
    \ncohomologyh{n}{G}{\field}
    &= \catstmod{G}(\Omega^n\field,\field) \\
    &\simeq \lim_{P \in \catop{\cattransport{G}{S}}}\catstmod{P}(\Omega^n\field,\field) \\
    &= \lim_{P \in \catop{\cattransport{G}{S}}}\ncohomologyh{n}{P}{\field}
  \end{align*}

  The Cartan-Eilenberg formula for the usual group cohomology
  $\cohomology{G}{\field}$ can be similarly recovered from the $2$-functor $\catder{-}$.

\end{example}

\subsection{Categorical invariants}

Another way to exploit \cref{theoremMackeyAs2Lim} is to try to factor
$\ffunctor{M}$ as a composite of two $2$-functors
\[
  \ffunctor{M} \colon \catop{(\catgpdf)} \xrightarrow{\tilde{\ffunctor{M}}} \bicatname{C} \xrightarrow{W} \catcat
\]
where $\tilde{\ffunctor{M}}$ is product preserving, $W \colon \bicatname{C} \to \catcat$
reflects bilimits, and to find a bilimit-preserving $2$-functor
\[
  H \colon \bicatname{C} \to \catname{D}
\]
to a $1$-category $\catname{D}$.
In this case, for any finite group $G$ with $p$-Sylow $S$:
\begin{equation}\label{eqInvLimit}
  \begin{split}
    H \circ \tilde{\ffunctor{M}}(G)
    & \simeq H(\llim_{\catop{\cattransportext{G}{S}}}\tilde{\ffunctor{M}} \circ U) \\
    & \simeq \llim_{\catop{\cattransportext{G}{S}}}H \circ \tilde{\ffunctor{M}} \circ U \\
    & \simeq \lim_{\catop{\cattransport{G}{S}}}H \circ \tilde{\ffunctor{M}} \circ U \qquad \text{(by \cref{propBilim1Cat})}
  \end{split}
\end{equation}

\begin{example}
  We sketch an example of application, taking $\ffunctor{M} = \catstmod{-}$. We can
  factor it through the $2$-category $\catadd_{\bullet}^{\text{gr}}$ of graded additive
  categories with a distinguished object, choosing the trivial module $\field$
  in each $\catstmod{G}$ and graduating by the Heller shift~$\Omega$. The forgetful functor $\catadd_{\bullet}^{\text{gr}} \to
  \catcat$ reflects bilimits. Then consider the well-defined bilimit-preserving
  2-functor to the $1$-category of rings:
  \[
    H \colon \left\{
      \begin{array}{lll}
        \catadd_{\bullet}^{\text{gr}} & \to & \textnormal{Ring} \\
        (C, \bullet_C) & \mapsto & C^*(\bullet_C,\bullet_C) \\
        (F : C \to D, u : F(\bullet_C) \xrightarrow{\sim} \bullet_D) & \mapsto & D(u^{-1},u)
                                                                     \circ F \\
        \alpha & \mapsto & \id
      \end{array}
    \right.
  \]
  Note that, for any finite group $G$:
  \[
    H \circ \tilde{\ffunctor{M}}(G) = \catstmod{G}^*(\field,\field) =
    \cohomologyh{G}{\field}
  \]
  Hence by \cref{eqInvLimit}, we obtain the same Cartan-Eilenberg formula as in \cref{exampleCEHom}:
  \begin{align*}
    \cohomologyh{G}{\field}
    &= H \circ \tilde{\ffunctor{M}}(G) \\
    &\simeq \lim_{\cattransport{G}{S}}H \circ \tilde{\ffunctor{M}} \circ U \\
    &= \lim_{P \in \cattransport{G}{S}}\cohomologyh{P}{\field}
  \end{align*}
\end{example}

It would be interesting to find other factorizations and bilimit-preserving $2$-functors $H$ with
value in a $1$-category for the various $p$-monadic Mackey 2-functors: this
would give possibly new ``Cartan-Eilenberg formulas''.

\printbibliography

\end{document}